\documentclass[11pt,openany]{book}
\usepackage{epsfig}

\usepackage{amssymb}

\newcommand{\pct}[1]{}
\newcommand{\Psfig}[1]{{\mbox{$\ \ $}}}

\begin{document}
\renewcommand{\thechapter}{\Roman{chapter}}

\thispagestyle{empty}

\
\vspace{0.4in}
 \begin{center}
 {\LARGE\bf KNOTS}\\
{\bf From combinatorics of knot diagrams to combinatorial topology
based on knots}
\end{center}

\vspace*{0.2in}

\centerline{Warszawa, November 30, 1984 -- Bethesda, March 3, 2007}
\vspace*{0.2in}

 \begin{center}
                      {\LARGE \bf J\'ozef H.~Przytycki}
\end{center}

\vspace*{0.2in}

\ \\
{\LARGE  LIST OF CHAPTERS}:\ \\
\ \\
{\LARGE \bf Chapter I: \ Preliminaries }\\
\ \\
{\LARGE \bf Chapter II:\ History of Knot Theory}\\
{\bf This e-print. Chapter II starts at page 3}\\
\ \\
{\LARGE \bf Chapter III:\ Conway type invariants }\\
\ \\
{\LARGE \bf Chapter IV:\  Goeritz and Seifert matrices}\\ \ \\
{\LARGE \bf Chapter V:\ Graphs and links}\\
{\bf e-print: http://arxiv.org/pdf/math.GT/0601227}\\
\ \\
{\LARGE \bf Chapter VI:\ Fox $n$-colorings, Rational moves, Lagrangian tangles
and Burnside groups}\ \\
\ \\
{\LARGE \bf Chapter VII:\ Symmetries of links}\ \\
\ \\
{\LARGE \bf Chapter VIII:\ Different links with the same
Jones type polynomials}\ \\
\ \\
{\LARGE \bf Chapter IX:\ Skein modules} \\
{\bf e-print: http://arxiv.org/pdf/math.GT/0602264}\\
\ \\
{\LARGE \bf Chapter X:\ Khovanov Homology: categorification of the Kauffman
bracket relation}\\
{\bf e-print: http://arxiv.org/pdf/math.GT/0512630 }\\ \ \\
{\LARGE \bf Appendix I.\ }\ \\
\ \\
{\LARGE \bf Appendix II.\ }\\ \ \\
{\LARGE \bf Appendix III.\ }\\
\

\ \\
{\LARGE \bf Introduction}\\
\ \\
This book is
about classical Knot Theory, that is, about
the position of a circle (a knot) or of a number of disjoint circles
(a link) in the space $R^3$ or in the sphere $S^3$.
We also venture into Knot Theory in general 3-dimensional
manifolds.

The book has its predecessor in Lecture Notes on Knot Theory,
which were published in Polish\footnote{The
Polish edition was prepared for the ``Knot Theory" mini-semester
at the Stefan Banach Center, Warsaw, Poland, July-August, 1995.}
in 1995  \cite{P-18}.
A rough translation of the Notes (by J.Wi\'sniewski) was
ready by the summer of 1995. It differed from the Polish edition
with the addition of
the full proof of Reidemeister's theorem. While I couldn't find
time to refine the translation and prepare the final manuscript,
I was adding new material and rewriting existing
chapters. In this way I created a new book based on the Polish
Lecture Notes
but expanded 3-fold.
Only the first part of Chapter III (formerly Chapter II),
on Conway's algebras is essentially unchanged from the Polish book
and is based on preprints \cite{P-1}.

As to the origin of the Lecture Notes, I was teaching an advanced course
in theory of 3-manifolds and Knot Theory at Warsaw University and it
was only natural to write down my talks (see Introduction to  (Polish)
Lecture Notes).
\\ \ \\
...\\
\ \  \ \ SEE Introduction before CHAPTER I.

\setcounter{chapter}{1}

\chapter{History of Knot Theory}\label{II}
 \vspace{0.2in}
{\LARGE \bf Abstract.}  Leibniz wrote in 1679: ``I consider that we need yet another
kind of analysis, $\ldots$
which deals directly with position." He called it
``geometry of position"(geometria situs). The first convincing example
of geometria situs was Euler's solution to the bridges of K\"onigsberg 
problem (1735). 
The first mathematical paper which mentions knots
was written by A.~T.~Vandermonde in 1771: ``Remarques sur les problemes 
de situation". We sketch in this chapter
the history of knot theory from Vandermonde to Jones stressing the
combinatorial aspect of the theory that is so visible in Jones type invariants.
In the first section we outline some older developments
 leading to modern knot theory.
\vspace*{0.2in}
\begin{quote}
{\em
``When Alexander reached Gordium,  he was seized with a longing to ascend to
the acropolis,   where the palace of Gordius and his son Midas was situated, 
and to see Gordius' waggon and the knot of the waggon's yoke$\ldots$. Over and
above this there was a legend about the waggon,  that anyone who untied the
knot of the yoke would rule Asia. The knot was of cornel bark,  and you could
not see where it began or ended. Alexander was unable to find how to untie
the knot but unwilling to leave it tied,  in case this caused a disturbance
among the masses; some say that he struck it with his sword,  cut the knot, 
and said it was now untied - but Aristobulus says that he took out the
pole-pin,  a bolt driven right through the pole,  holding the knot together, 
and so removed the yoke from the pole. I cannot say with confidence what
Alexander actually did about this knot,   but he and his suite certainly
left the waggon with the impression that the oracle about the undoing
of the knot had been fulfilled,  and in fact that night there was thunder
and lightning,   a further sign  from heaven; so Alexander in thanksgiving
offered sacrifice next day to whatever gods had shown the signs and
the way to undo the knot."
}
{\footnotesize \ \ \ 
[Lucius Flavius Arrianus, Anabasis Alexandri, Book II, c.150 A.D.,\cite{Arr}]}
\end{quote}
Similar account, clearly based on the same primary sources is given by
Plutarch of Chareonera (c. 46 - 122 A.D.). He writes in his 
``Lives" \cite{Plu} (page 271):
\begin{quote}
{\em
``Next he marched into Pisidia where he subdued any resistance which
he encountered, and then made himself master of Phyrgia. When he captured 
Gordium [in March 333 B.C.] which is reputed to have been the home of
the ancient king Midas, he saw the celebrated chariot which was
fastened to its yoke by the bark of the cornel-tree, and heard the
legend which was believed by all barbarians, that the fates had
decreed that the man who untied the knot was destined to become 
the ruler of the whole world. According to most writers the fastenings
 were so elaborately intertwined and coiled upon one another that their 
ends were hidden: in consequence Alexander did not know what to do,
and in the end loosened the knot by cutting through it with his
sword, whereupon the many ends sprang into view. But according to
Aristobulus he unfastened it quite easily by removing the pin which
secured the yoke to the pole of the chariot, and then pulling out
the yoke itself." }

\end{quote}

\newpage

   In this chapter we present the history of ideas which lead up
to the development of modern knot theory. We are more 
detailed when pre-XX
century history is reported. With more recent times we are more selective,  
stressing developments related to Jones type invariants of links.
Additional historical information on specific topics of Knot Theory 
is  given in other chapters of the book\footnote{There are 
books which treat the history of topics related to knot theory
\cite{B-L-W,Ch-M,Crowe,Die,T-G}. 
J.Stillwell's textbook \cite{Stil} contains very interesting historical
digressions.}.

Knots  have fascinated people from the dawn of the human history. 
We can wonder what caused a merchant living about 1700 BC. in Anatolia
and exchanging goods with Mesopotamians, to choose braids and knots as
his seal sign; Fig.1.1. We can guess however that 
stamps, cylinders and seals with knots
and links as their motifs appeared before proper writing was invented 
about 3500 BC.
\ \\
\ \\
\centerline{\psfig{figure=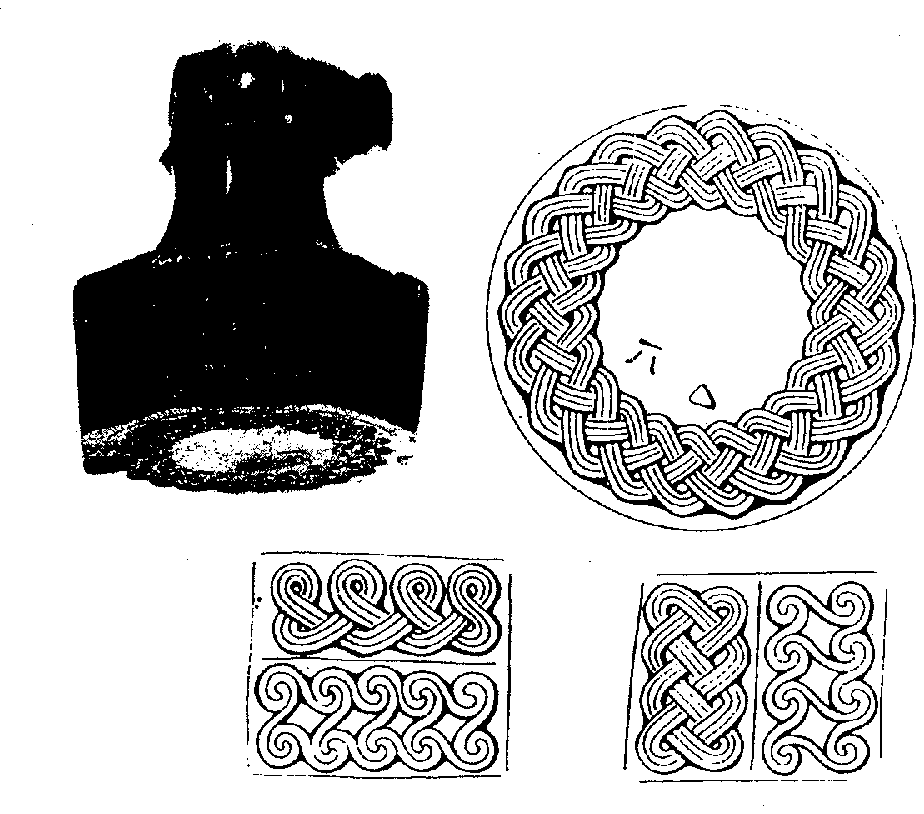,height=8.4cm}}
\begin{center} 
                Figure 1.1 Stamp seal, about 1700 BC (the British Museum).
\end{center} 
{\footnotesize {\it On the octagonal base [of hammer-handled haematite seal]
 are patterns surrounding a hieroglyphic inscription (largely erased).
Four of the sides are blank and the other four are engraved with elaborate
patterns typical of the period (and also popular in Syria) alternating
with cult scenes...}(\cite{Col}, p.93). }

I am unaware of any pre-3500 BC examples but I will describe two finds 
from the third millennium BC.

The oldest examples  outside Mesopotamia, that I am aware of,
 are from  the pre-Hellenic Greece.
Excavations at Lerna
by the American School of Classical Studies under the direction of
Professor J.~L. Caskey (1952-1958) discovered two rich
deposits of clay seal-impressions.  The second deposit dated
from about 2200 BC
contains several impressions of knots and links\footnote{The early Bronze
Age in Greece is divided,
as in Crete and the Cyclades, into three phases. The second phase lasted
from 2500 to 2200 BC, and was marked by a considerable increase in
prosperity. There were palaces at Lerna, and Tiryns, and probably elsewhere,
in contact with the Second City of Troy. The end
of this phase (in the Peloponnese)
was brought about by invasion and mass burnings.
The invaders are thought to be the first speakers of
the Greek language to arrive in Greece.} \cite{Hig,Hea,Wie}
 (see Fig.1.2). \\
\ \\

\centerline{\psfig{figure=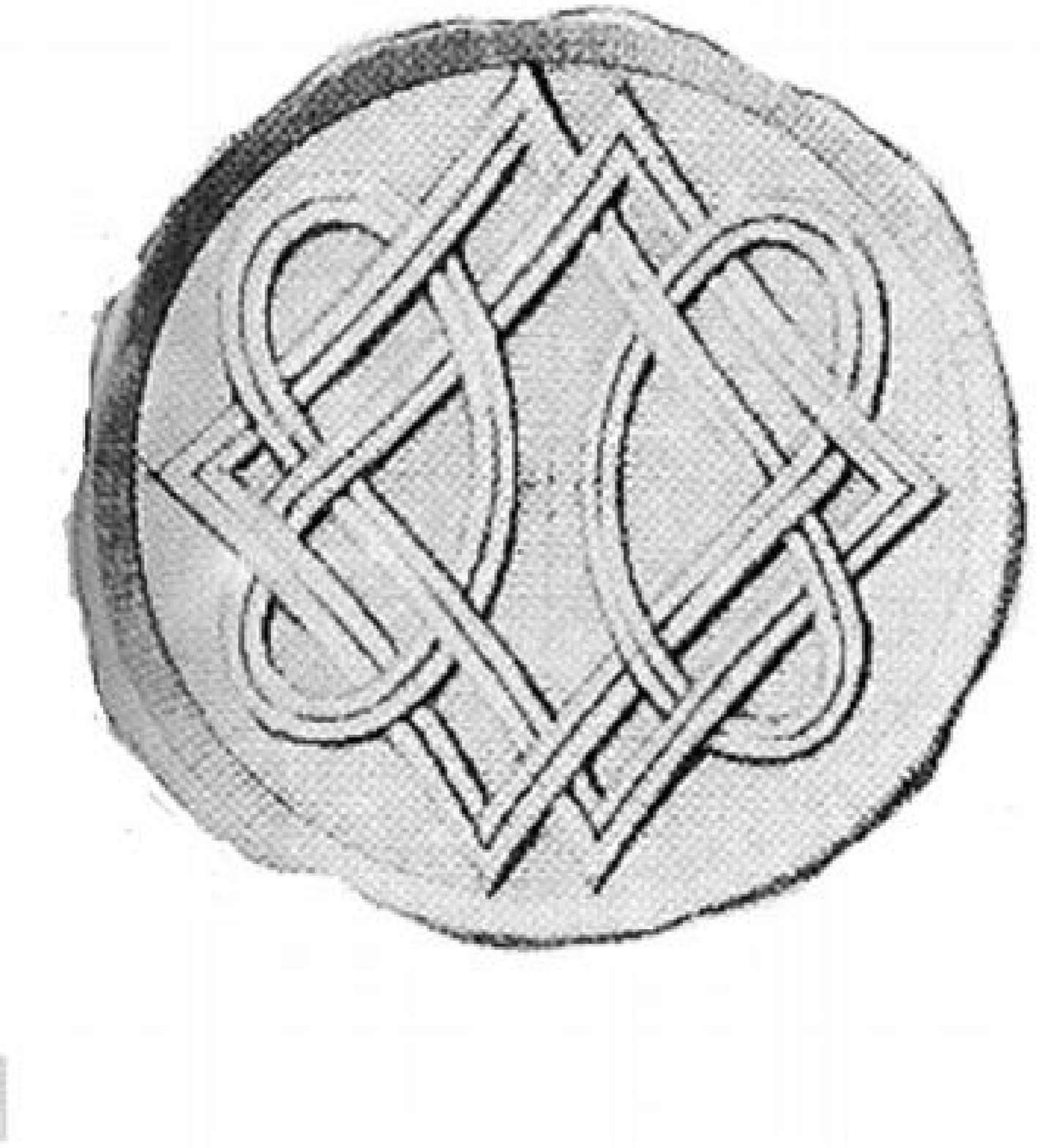,height=10.5cm}}
\centerline{Fig. 1.2; \ A seal-impression from the House of
the Tiles in Lerna \cite{Hig}.}
\ \\

Even older example of cylinder seal impression (c. 2600-2500 B.C.)
 from Ur, Mesopotamia is described  
in the book Innana by Diane Wolkstein and 
Samuel Noah Kramer \cite{Wo-Kr} (page 7); Figure 1.3,  
illustrating the text:

``Then a serpent who could not be charmed
made its nest in the roots of the tree."
\ \\
\ \\
\centerline{\ \ \ \ \ \ \ \ \ \ \ \ \ 
\psfig{figure=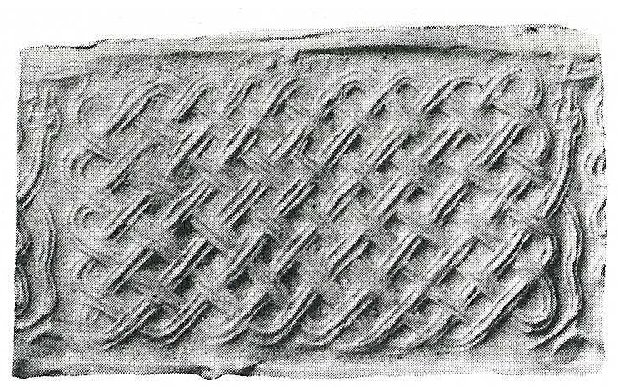,height=8.5cm}}
\centerline{Fig. 1.3; \ {\it Snake with Interlacing Coil}}
{\footnotesize {\it 
Cylinder seal. Ur, Mesopotamia. The Royal Cemetery, Early Dynastic
period, c. 2600-2500 B.C. Lapis lazuli. Iraq Museum. Photograph
courtesy of the British Museum, UI 9080, 
\cite{Wo-Kr}\footnote{On the pages 179-180 they comment:\ 
   The majority of the pictorial surface is covered with the inter-
twined coils of a serpent, forming a lattice pattern. To the right its tail
appears below the coils and its head above, with a bird perched upon it.
   Two snakes intertwined rather than one are shown on earlier
representations of this motif. Snakes twist themselves together in this
fashion when mating, suggesting this symbol's association with fertility.}}}

\section{From Heraclas to D\"urer}\label{Section II.2}
 
   It is tempting to look for the origin of knot theory in Ancient
Greek mathematics (if not earlier). There is some justification to do so:
a Greek physician named Heraklas, who lived during the first century A.D.
and who was likely a pupil or associate of Heliodorus, wrote an essay
on surgeon's slings\footnote{Heliodorus, who lived at the time 
of Trajan (Roman Emperor 98--117 A.D.), also mentions in 
his work knots and loops \cite{Sar-1}}\footnote{Hippocrates of 
Cos (c.460 - 375 B.C.) in his collection of notes:
In the surgery; De officina medici; Cat' i\={e}treion, deals with
bandaging. Thessalos, Hippocrates' son, has been named also as the author.
 A commentary on the Hippocratic treatise on {\em Joints} was written by
 Apollonios of Citon (in Cypros), who flourished in Alexandria in the
first half of the first century B.C. That commentary has obtained a
great importance because of an accident in its transmition. A manuscript
 of it in Florence (Codex Laurentianus) is a Byzantine copy of the ninth
century, including surgical illustrations (for example, with reference to
 reposition methods), which might go back to the time of Apollonios and
even Hippocrates. Iconographic tradition of this kind are very rare,
because the copying of figures was far more difficult than the writing
of the text and was often abandoned \cite{Sar-1}.
The story of the illustrations to Apollonios' commentary is described in
\cite{Sar-2}.}.
Heraklas explains,  giving step-by-step instructions, 
eighteen ways to tie orthopedic slings. His work survived because
Oribasius of Pergamum (ca. 325-400; physician of the emperor
 Julian the Apostate) included it toward the end of the fourth century
in his ``Medical Collections". The oldest extant manuscript of ``Medical
Collections" was made in the tenth century by the Byzantine physician
Nicetas. The Codex of Nicetas was brought to Italy in the fifteenth century by 
an eminent Greek scholar,   J.~Lascaris,  a refugee from Constantinople.
Heraklas' part of the Codex of Nicetas has no illustrations,  
and around 1500 an anonymous artist depicted Heraklas' knots in one of
the Greek manuscripts of Oribasus ``Medical Collections" (in Figure 2.1 we
reproduce the drawing of the third Heraklas knot together with its 
original,   Heraklas',  description). Vidus Vidius (1500-1569),  
a Florentine who became physician to Francis I (king of France, 1515-1547)
 and professor of medicine in the Coll\`ege
de France,   translated the Codex of Nicetas into Latin; it contains also
drawings of Heraklas' surgeon's slings by the Italian painter,  
sculptor and architect Francesco Primaticcio (1504-1570); \cite{Da,Ra}.
\vspace*{1in}
\centerline{\psfig{figure=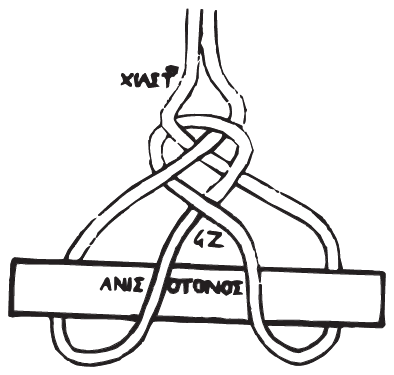,,height=7.6cm}} 
 \centerline{  Figure 2.1; The crossed noose}
{\footnotesize ``For the tying the crossed noose,  a cord,
folded double,  is procured,
and the ends of the cord are held in the left hand,   and the loop is
held in the right hand. Then the loop is twisted so that the slack
parts of the cord crossed. Hence the noose is called crossed. After
the slack parts of the cord have been crossed,  the loop is placed
on the crossing,   and the lower slack part of the cord is pulled up
through the middle of the loop. Thus the knot of the noose is in the
middle,   with a loop on one side and two ends on the other. This
likewise,   in function,  is a noose of unequal tension";
\cite{Da}}.

   Heraklas' essay should be taken seriously as far as knot theory is
concerned even if it is not knot theory proper but rather its application.
The story of the survival of Heraklas' work; and efforts to reconstruct
his knots in Renaissance is typical of all science disciplines and
efforts to recover lost Greek books provided the important engine for
development of modern science. This was true in Mathematics as well:\
the beginning of modern calculus in XVII century can be traced
to efforts of reconstructing lost books of Archimedes and other ancient
Greek mathematicians. It was only the work of Newton and Leibniz which
went much farther than their Greek predecessors.

There are two enlighting examples of great Renaissance artists 
interest in knots:\ 
Engravings by Leonardo da Vinci\footnote{Giorgio Vasari writes
in \cite{Va}: ``[Leonardo da Vinci] spent much time in making a regular
design of a series of knots so that the cord may be traced
from one end to the other, the whole filling a round space.
There is a fine engraving of this most difficult design,
and in the middle are the words: Leonardus Vinci Academia."} (1452-1519)
 \cite{Mac}\ 
\ \\
\ \\
\centerline{\psfig{figure=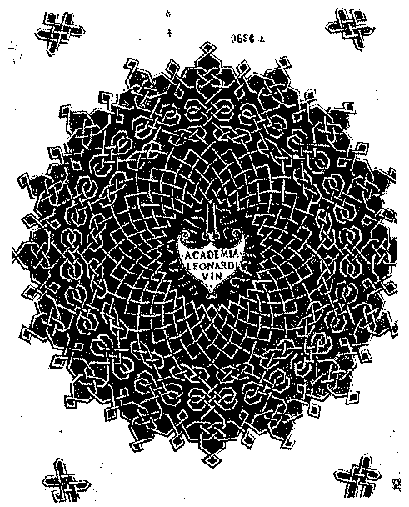,height=10.5cm}}
\centerline{Fig. 2.2; \ A knot by Leonardo \cite{Mac}; c. 1496}
  and woodcuts by Albrecht
D\"urer\footnote{``Another great artist with whose works D\"urer now became
acquainted was Leonardo da Vinci. It does not seem likely that the
two artists ever met, but he may have been brought into relation
with him through Luca Pacioli, the author of the book De Divina
Proportione, which appeared at Venice in 1509, and an intimate friend
of the great Leonardo. D\"urer would naturally be deeply interested in
the proportion theories of Leonardo and Pacioli. He was certainly
acquainted with some engravings of Leonardo's school, representing a
curious circle of concentric scrollwork on a black ground, one of them
entitled Accademia Leonardi Vinci; for he himself executed six woodcuts
in imitation, the Six Knots, as he calls them himself. D\"urer was
amused by and interested in all scientific or mathematical problems..."
From: http://www.cwru.edu/edocs/7/258.pdf, compare
\cite{Dur-2}.} (1471-1528) \cite{Dur-1,Ha}, Fig.2.3.\\
\ \\
\centerline{\psfig{figure=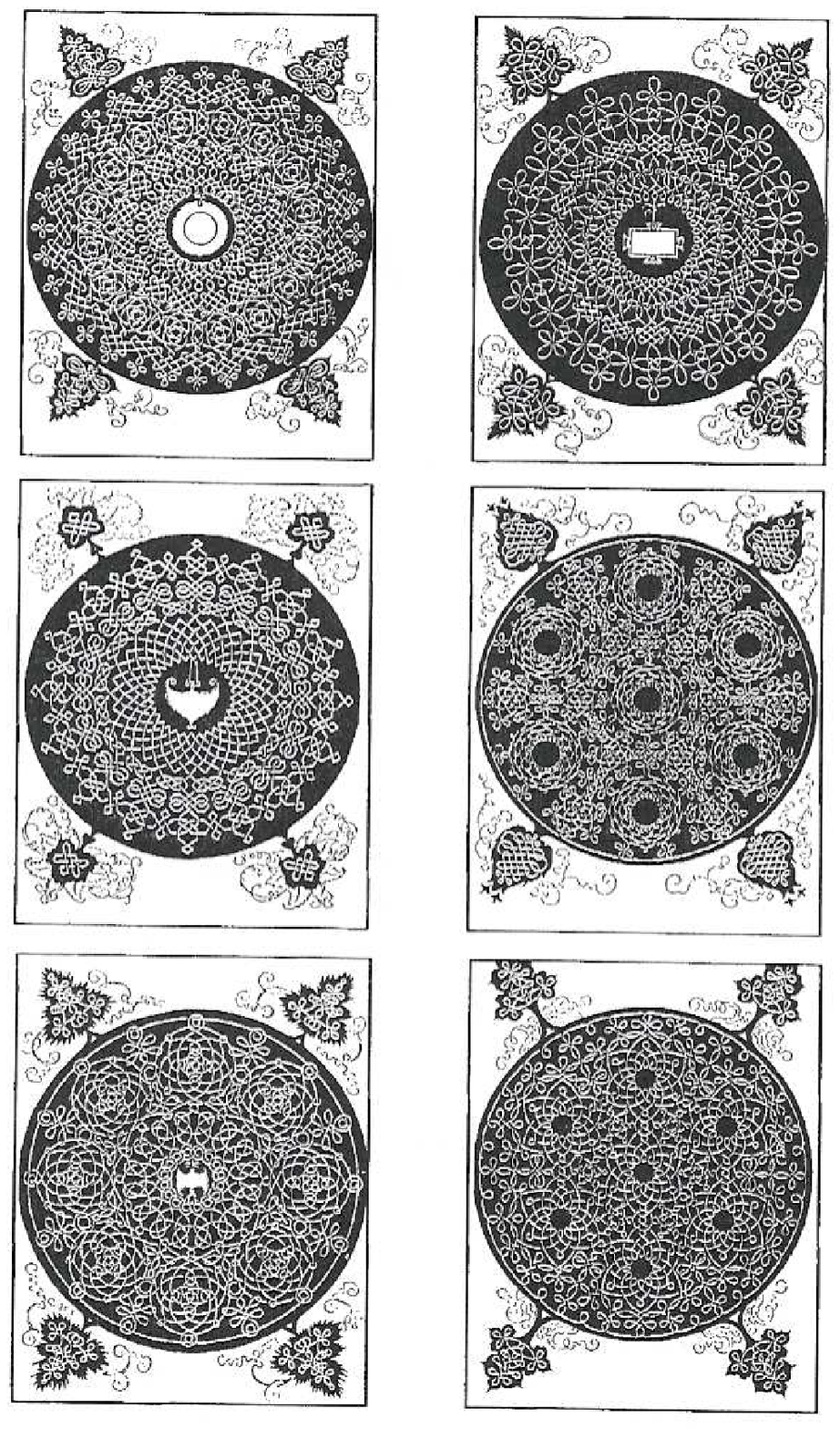,height=18.5cm}}
\centerline{Fig. 2.3; \ Six knots by D\"urer \cite{Kur}; c. 1505-1507}
\ \\

\section{Dawn of Knot Theory}\label{Section II.3}

We would argue, that modern knot theory has its roots with 
 Gottfried Wilhelm Leibniz (1646-1716) speculation that aside from calculus 
and analytical geometry there should exist a ``geometry of position"  
(geometria situs) which deals with relations depending on position 
alone (ignoring magnitudes). In a letter to Christian Huygens (1629-1695), 
written in 1679 \cite{Lei}, he declared:\\
``I am not content with algebra,  in that it yields neither the shortest 
proofs nor the most beautiful constructions of geometry. 
Consequently, in view of this, I consider that we need yet 
another kind of analysis, geometric or linear, which deals directly 
with position, as algebra deals with magnitude".

   I do not know whether Leibniz had any convincing example of a problem
belonging to the geometry of position. According to \cite{Kli}:\\ 
``As far back as 1679 Leibniz, in his Characteristica Geometrica, 
tried to formulate basic geometric properties of geometrical figures,  
to use special symbols to represent them,  and to combine these 
properties under operations so as to produce others. 
He called this study analysis situs or geometria
situs...  To the extent that he was at all clear, Leibniz envisioned what
we now call combinatorial topology".

   The first convincing example of geometria situs was studied by 
Leonard Euler (1707-1783). This concerns the bridges on the river Pregel 
at K\"onigsberg (then in East Prussia)\footnote{Euler never 
visited K\"onigsberg. 
He was informed about the puzzle of bridges of K\"onigsberg 
(and possible relation to Leibniz geometria situs) by 
future Major of Danzig (Gda\'nsk) Carl Leonhard 
Gottlieb Ehler (1685-1753); the first letter from Ehler to Euler is 
dated April 8, 1735.
Ehler in turn was acting on behalf of Danzig mathematician 
Heinrich Kuhn (1690-1769) \cite{H-W}. Kuhn was born in 
K\"onigsberg, he studied at the Pedagogicum there, 
... in 1733 he settled in Danzig. One should add that Kuhn 
was the first person to suggest geometric interpretation of complex 
numbers \cite{Jan}.}. 
 Euler solved (and generalized) the bridges of K\"onigsberg problem
and on August 26, 1735 presented his solution to the Russian Academy at
St. Petersburg (it was published in 1736), \cite{Eu}. 
With the Euler paper, graph theory and topology were born. 
Euler started his paper by remarking:\\
``The branch of geometry that deals with magnitudes has been zealously studied
throughout the past,  but there is another branch that has been almost 
unknown up to now; Leibniz spoke of it first,  calling it the 
``geometry of position" (geometria situs). 
This branch of geometry deals with relations dependent
on position; it does not take magnitudes into considerations, nor does it
involve calculation with quantities. But as yet no satisfactory definition
has been given of the problems that belong to this geometry of position
or of the method to be used in solving them".


\begin{center}
\centerline{\psfig{figure=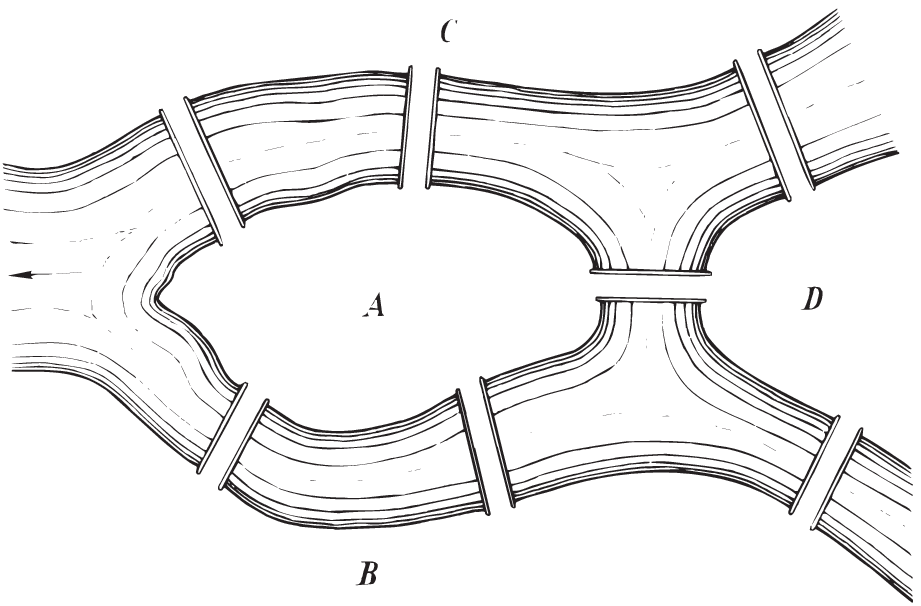}}
Figure 3.1; Bridges of K\"onigsberg
\end{center}

   For the birth of knot theory one had to wait another 35 years. In 1771
Alexandre-Theophile Vandermonde (1735-1796) wrote the paper: 
{\it Remarques sur les probl\`emes de situation} (Remarks on
problems of positions) where he specifically places braids and knots as
a subject of the geometry of position \cite{Va}. 
In the first paragraph of the paper
Vandermonde wrote: 

{\it Whatever the twists and turns of a system of threads
in space, one can always obtain an expression for the calculation of its
dimensions, but this expression will be of little use in practice.
The craftsman who fashions a braid, a net, or some knots will be concerned, 
not with questions of measurement, but with those of position: what he sees
there is the manner in which the threads are interlaced.} 

\centerline{\psfig{figure=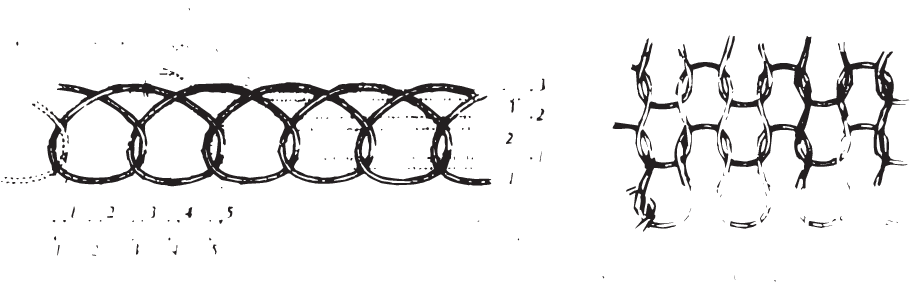}}
\begin{center}
Figure 3.2; Knots of Vandermonde
\end{center}

   In our search for the origin of knot theory, we arrive next at
Carl Friedrich Gauss (1777-1855). According to  \cite{Stac,Dun}
:

``One of the oldest notes by Gauss to be found among his papers is a sheet
of paper with the date 1794. It bears the heading ``A collection of knots"
and contains thirteen neatly sketched views of knots with English names
written beside them... With it are two additional pieces of paper with
sketches of knots. One is dated 1819; the other is much later,  ...".
\footnote{According to \cite{Gr-2}, the first English sailing book with
pictures of knots appeared in 1769 \cite{Falc}. 
}

In July of 1995 I finally visited the old library in G\"ottingen,
I looked at knots from 1794 - in fact not all of them are drawn - some
only described; see Fig. 3.3 for one of the drawings\footnote{First 
eight drawings are reproduced in the preface to \cite{T-G}.}. 
\ \\
\ \\
\centerline{\psfig{figure=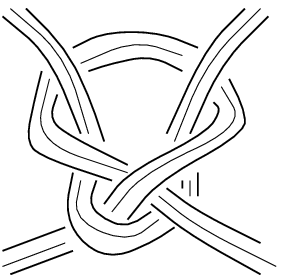,height=5.1cm}} 
\begin{center} 
     Figure 3.3; {\it Meshing knot}, 10'th knot of Gauss from 1794. 
\end{center} 
There are other fascinating drawings in Gauss' notebooks.
For example, the drawing of a braid with complex coordinate 
description at each height (Figure 3.4; compare \cite{Ep-1, P-21}), 
 and the note that it is a good method of coding a knotting.
It is difficult to date the drawing; one can say for sure that it was done
between 1814 and 1830, I would guess closer to 1814\footnote{As a curiosity 
one can add that of one of the notebooks (Handb. 3) in which Gauss had 
also drawn some knot diagrams has braids motives on its cover.}.
\ \\ 
\ \\
\centerline{\psfig{figure=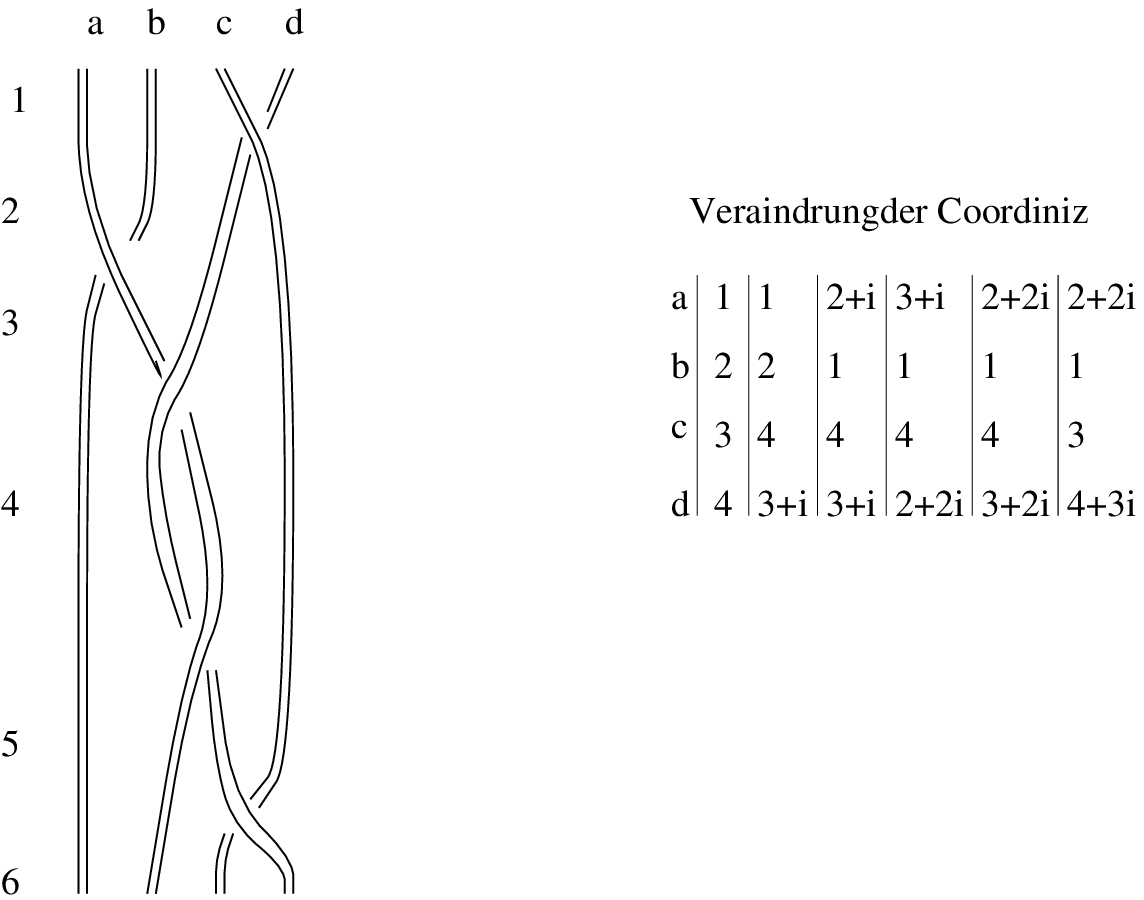}} 
\begin{center} Figure 3.4\footnote{Gauss coordinates are not always 
consistent; most of the time $i$ is pointing downward but there are exceptions.}
 \end{center} 
\begin{center} 
{\it It is a good method of coding a knotting} (from 
a Gauss' notebook (Handb.7)).
\end{center}

  There is also the mysterious ``framed tangle",
 see Fig.3.5 \cite{Ga-1,P-29} whose interpretation is not yet 
convincingly given.

\centerline{\psfig{figure=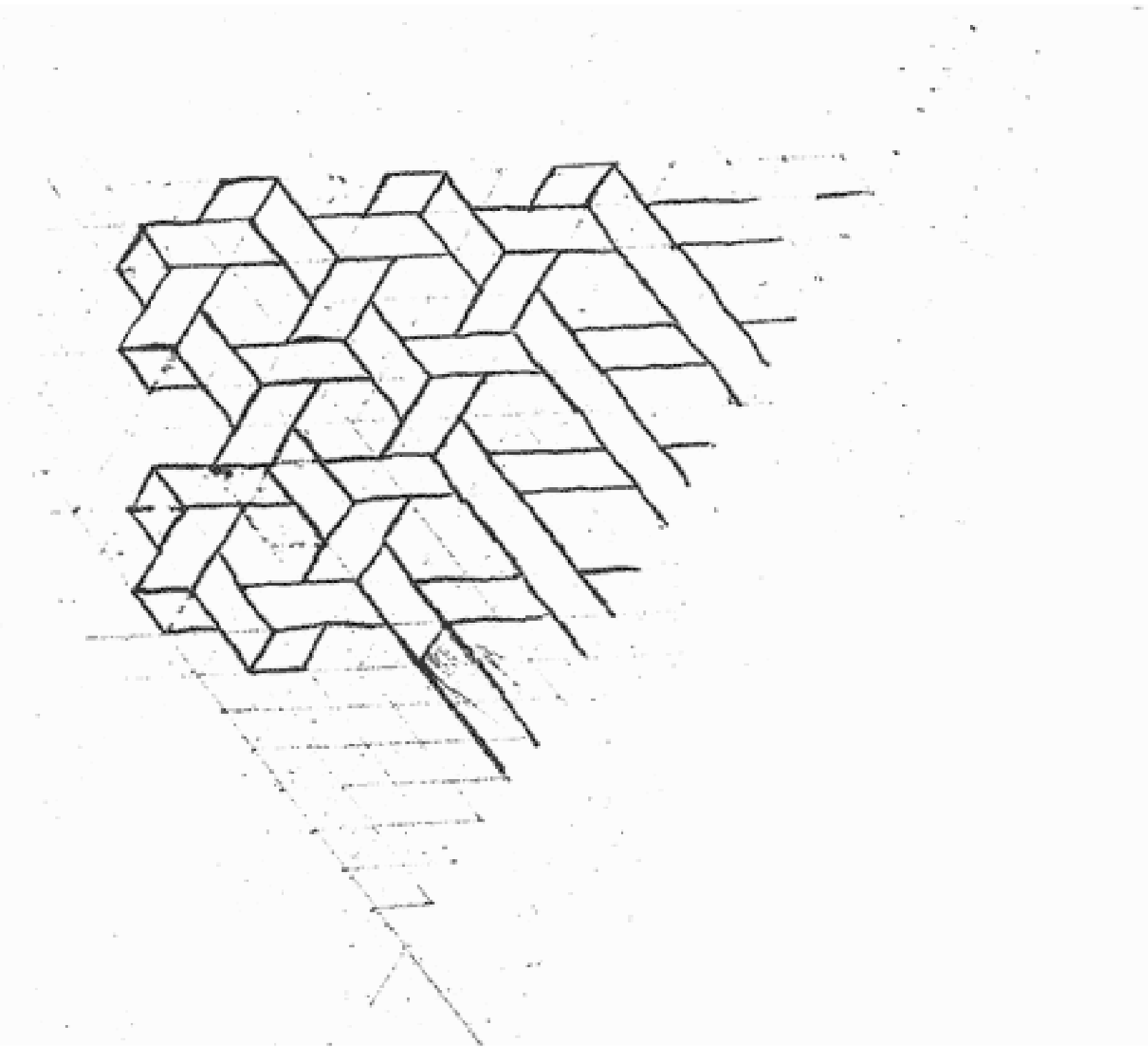,height=10.6cm}}
\centerline{Fig. 3.5; \ Framed tangle from Gauss' notebook \cite{Ga-1}}

   In his note (Jan. 22 1833) Gauss introduces the linking number of two 
knots\footnote{His method is analytical - the Gauss integral; 
in modern language Gauss integral computes analytically the degree of 
the map from a torus parameterizing a 2 component link to the unit 2-sphere.}. 
Gauss' note presents the first deep incursion into knot theory; 
it establishes that the following two links
are substantially different:\\ \parbox{1.2cm}{\psfig{figure=T2.eps}}\ , 
\ \parbox{0.9cm}{\psfig{figure=+Hopfmaly.eps}}.  \ \ 
Gauss' analytical method has
been recently revitalized by Witten's approach to knot theory \cite{Wit}.

James Clerk Maxwell (1831-1879), in his fundamental book of 1873 
``A treatise on electricity \& magnetism" \cite{Max} writes
\footnote{It was only six years after 
Gauss note was first published in his collected works in 1867.}
: 
``It was the discovery by Gauss of
this very integral, expressing the work done on a magnetic pole while
describing a closed curve in presence of a closed electric current,
and indicating the geometrical connection between the two closed curves,
that led him to lament the small progress made in the Geometry
of Position since the time of Leibnitz, Euler and Vandermonde. We
have now, however, some progress to report chiefly due to Riemann,
Helmholtz and Listing."\footnote{Gauss wrote in 1833, in the same note in
which he introduced the linking number: 
``On the geometry of position, which Leibniz initiated and to which only two
geometers,  Euler and Vandermonde,  have given a feeble glance,  we know and
possess, after a century and a half, very little more than nothing."}
Maxwell 
goes on to describe two closed curves which cannot be separated
but for which the value of the  Gauss integral is equal to zero; Fig.3.6.

\centerline{\psfig{figure=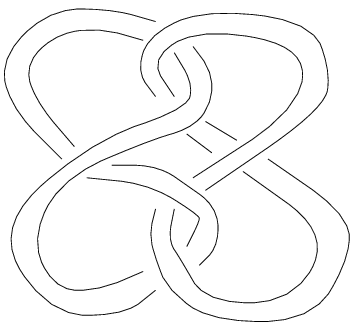}}
\begin{center}
Figure 3.6;\ The link of Maxwell.
\end{center}

In 1876, O.Boeddicker observed that, in a certain sense, the linking
number is the number of the crossing points of the second curve with
a surface bounded by the first curve \cite{Boe-1,Boe-2,Bog}.\
Hermann Karl Brunn\footnote{Born August 1, 1862 Rome (Italy), 
died Sept. 20, 1939 (Munchen, Germany)\cite{Bl}.}
 \cite{Br} observed in 1892 
that the linking number 
of a two-component link, considered by Gauss, can be read from 
a diagram of the link\footnote{It is also noted by Tait 
in 1877 (\cite{Ta}, page 308).}. 
If the link has components $K_1$ and $K_2$, we
consider any diagram of the link and count each point
at which $K_1$ crosses under $K_2$ as $+1$ for 
\parbox{0.6cm}{\psfig{figure=L+maly.eps}}\
and $-1$ for\ \parbox{0.6cm}{\psfig{figure=L-maly.eps}}\ . 
The sum of these, over all crossings
of $K_1$ under $K_2$, is the Gauss linking number.
 
   The year of 1847 was very important for the knot theory 
(graph theory and topology as well). 
On one hand,  Gustav Robert Kirchhoff (1824-1887)
published his fundamental paper on electrical circuits \cite{Kir}. 
It has deep connections with knot theory,  
however the relations were discovered only about
a hundred years later (e.g. the Kirchhoff complexity of a circuit corresponds
to the determinant of the knot or link determined by the circuit).
On the other hand, Johann Benedict Listing (1808-1882), a student of Gauss, 
published his monograph (Vorstudien zur Topologie, \cite{Lis}). 
A considerable part
of the monograph is devoted to knots. Even earlier, on April 1, 1836,
 Listing wrote a letter from Catania to "Herr Muller", his former school 
teacher\footnote{Johann Heinrich M\"uller (1787-1844) was the mathematics 
and astronomy master at {\it Musterschule} which Listing entered in 1816 
\cite{Bre-2}.}, with the heading 
"Topology", concerning ... 
(2) winding paths of knots; and (3) paths in a lattice 
\cite{Bre-1,Bre-2,Stac}. 
Listing stated in particular that the right handed trefoil 
knot ( \parbox{1.2cm}{\psfig{figure=+trefoilmaly.eps}}\ ) 
 and the left handed trefoil knot 
(\parbox{1.2cm}{\psfig{figure=min-trefoilmaly.eps}}) 
are not equivalent. Later Listing showed that
the figure eight knot \ \parbox{1.2cm}{\psfig{figure=+figureeightmaly.eps}}\  
and its mirror image \ 
\parbox{1.2cm}{\psfig{figure=min-figureeightmaly.eps}}\    are equivalent
(we say that the figure eight knot, also called the Listing knot, is
amphicheiral)\footnote{This was observed in the note dated March 18, 1849
\cite{Lit-1}}. 
Listing indebtness to Gauss is nicely described in the 
introduction to \cite{Lis}:\ \ 
``Stimulated by by the greatest geometer of our days, 
who had been repeatedly turning my attention to the significance of this
subject, during long time I did various attempts to analyze different cases
related to the subject, given by natural sciences and their applications.
And if now, when these reflections do not have a right yet to claim
rigorous scientific form and method, I let myself to publish them as
preliminary sketch of the new science, then I do this with the intention
to turn attention to significance and potential of it, with help of collected
here important information, examples and materials.
   I hope you let me use the name ``Topology" for this kind of studies
of spatial images, rather than suggested by Leibniz name "geometria situs",
reminding of notion of ``measure" and "measurement", playing here absolutely
subordinate role and confiding with ``g\'eom\'etrie de position" which
is established for a different kind of geometrical studies.
Therefore,  by {\bf Topology} we will mean the study of modal relations of
spatial images, or of laws of connectedness, mutual disposition and traces
of points, lines surfaces, bodies and their parts or their unions in space,
independently of relations of measures and quantities. By means of the
notion "trace", which is very close to the notion of movement, topology
is related to mechanics, similarly as it is related to geometry.
Of course, velocity, as well as mass, momentum, powers and moments of
movement from the quantity point of view are not taken into consideration.
Instead we consider only modal relations between moving or moved in
space images. In order to reach the level of exact science, topology will
have to translate facts of spatial contemplation into easier notion
which, using corresponding symbols analogous to mathematical ones, we
will be able to do corresponding operations following some simple rules."
{\small (Translated by M.Sokolov)}.

As we mentioned before, Maxwell, in his study of electricity and magnetism,  
had some thoughts on knots and links (in particular motivated by the 
freshly published Gauss' collected works). The origin of modern knot
theory should be associated with four physicists: Hermann Von Helmholtz 
(1821-1894), William Thomson (Lord Kelvin) (1824-1907), Maxwell 
and Peter Guthrie Tait (1831-1901). 
We can quote after Tait's assistant in Edinburgh and later biographer,
C.G.Knott \cite{Kno}:

{\it Tait was greatly impressed with Helmholtz's famous paper on 
vortex motion [\cite{Helm}; 1858]... Early in 1867 he devised 
a simple but effective method of producing vortex smoke rings;
and it was when viewing the behaviour of these in Tait's Class Room that
Thomson was led to the conception of the vortex atom. In his first paper
to the Royal Society of Edinburgh on February 18, 1867 [\cite{Thoms}], 
Sir William Thomson refers... to the genesis of the conception.} 
In turn Thomson's theory was Tait's motivation to understand the structure
of knots. In Tait's words: {\it I was led to the consideration of the form of
knots by Sir W.~Thomson's Theory of Vortex Atoms, and consequently the point
of view which,  at least at first, I adopted was that of classifying knots
by the number of their crossings... The enormous number of lines in the spectra
of certain elementary substances show that, if Thomson's suggestion
be correct, the form of the corresponding vortex atoms cannot be regarded
as very simple[\cite{Ta}].}

\begin{center}
\centerline{\psfig{figure=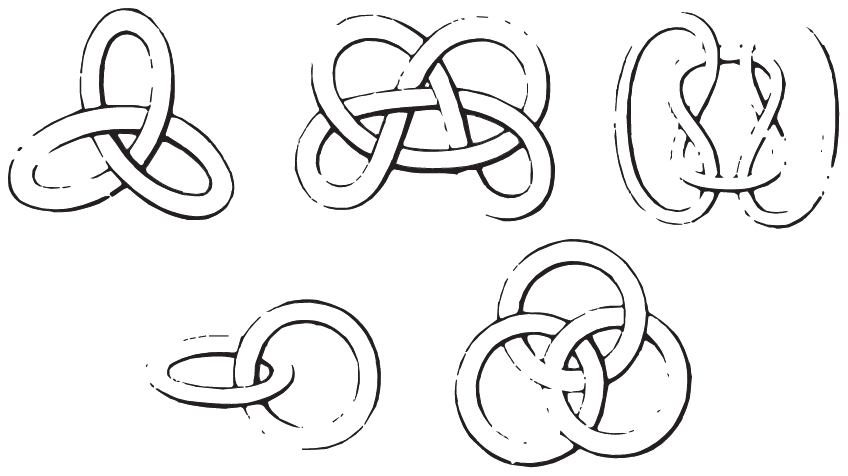}}
Figure 3.7. Knots and links of William Thomson (Kelvin) from 1867.
\end{center}
There is an interesting letter from Maxwell to Tait dated Nov. 13, 1867,
which shows that Tait was sharing his ideas of knots with his friend
\cite{Kno,Lom}. In one of his rhymes Maxwell wrote (clearly referring to Tait):
\\
\ \\
{\it Clear your coil of kinkings\\
\ \ \ \ Into perfect plaiting,\\
Locking loops and linkings\\
\ \ \ \ Interpenetrating.}\\
\cite{Kno}

Tait describes his work on knots in the following words (\cite{Ta},1877):
{\it When I commenced my investigations I was altogether unaware that
anything had been written (from a scientific point of view) about knots.
No one in Section A at the British Association of 1876, when I read a little
paper[\cite{Ta-0}] on the subject, could give me any reference; and it was not
till after I had sent my second paper to this Society that I obtained, 
in consequence of a hint from Professor Clerk-Maxwell, a copy of the very 
remarkable Essay by Listing\footnote
{In 1883 Tait wrote in {\it Nature} obituary after Listing death \cite{Ta-1}:
{\it One of the few remaining links that still continued to connect our
time with that in which Gauss had made G\"ottingen one of the chief
intellectual centres of the civilised world has just be broken by the death
of Listing... This paper [{\it Vorstudien zur Topologie}], which is
throughout elementary, deserves careful translation into English...}.
After more than a hundred years the paper is not translated 
(only Tait summary exists \cite{Ta-2}) and one should repeat
Tait appeal again: The paper very much deserves translation. One can add
that in 1932 the paper was translated into Russian.}, 
{\it Vorstudien zur Topologie}[\cite{Lis}],
of which (so far as it bears upon my present subject) I have given a full 
abstract in the Proceedings of the Society for Feb. 3, 1877. 
Here, as was to be expected, I found many of my
results anticipated, but I also obtained one or two  hints which,
though of the briefest, have since been very useful to me. Listing does not
enter upon the determination of the number of distinct form of knots with
a given number of intersections, in fact he gives only a very few forms
as examples, and they are curiously enough confined to three, five and
seven crossings only; but he makes several very suggestive remarks about
the representation of a particular class of ``reduced" knots... This work
of Listing's and an acute remark made by Gauss (which with some comments
on it by Clerk-Maxwell, will be referred to later), seem to be all of any
consequence that has been as yet written on the subject.} Tait's paper was
revised May 11, 1877; he finishes the paper as follows:
{\it After the papers, of which the foregoing is a digest, had been read, 
I obtained from Professor Listing\footnote{Library of the University
of California has a copy of {\it Vorstudien zur Topologie} which
Listing sent to Tait with the dedication.} and Klein a few references to the
literature of the subject of knots. It is very scanty, and has scarcely
any bearing upon the main question which I have treated above. Considering
that Listing's Essay was published thirty years ago, and that it seems to 
be pretty well known in Germany, this is a curious fact. From Listing's
letter (Proc.~R.S.E.. 1877, p.316), it is clear that he has published
only a small part of the results of his investigations. Klein himself
\cite{Klein} has made the very singular discovery that in space of four
dimensions there cannot be knots.\footnote{Klein observation was noticed 
in non-mathematical circles and it became part of popular culture. 
For example, the American magician and medium 
Henry Slade was performing ``magic tricks" claiming that he 
solves knots in fourth dimension. He was taken seriously by a German 
astrophysicist J.K.F.Zoellner who had with him a number of seances 
in 1877 and 1878.}

The value of Gauss's integral has been discussed at considerable length
by Boeddicker ... in an Inaugural Dissertation, with the title
{\em Beitrag zur Theorie des Winkels}, G\"ottingen, 1876.

An inaugural Dissertation by Weith, {\it Topologishe Untersuchung der
Kurven-Verschlingung}, Z\"urich, 1876 \cite{Weith}, is professedly based on 
Listing's Essay. It contains a proof that there is an infinite number 
of different forms of knots!\footnote{In fact it was proven only 20-30 years
later and depended on the fundamental work of Poincar\'e on foundation of
algebraic topology.} The author points out what he (erroneously) supposes
to be mistakes in Listing's Essay; and, in consequence, gives as something
quite new an illustration of the obvious fact that there can be irreducible
knots in which the crossing are not alternately over and under
\footnote{It was proven only in 1930 by Bankwitz \cite{Ban}, using
the determinant of a knot.}. The rest
of this paper is devoted to the relations of knots to Riemann's surfaces.}

   Tait, in collaboration with Reverend Thomas Penyngton Kirkman (1806-1895), 
and independently Charles Newton Little\footnote{Born Madura India, May 19, 
1858. A.B., Nebraska 1879, A.M. 1884; Ph.D, Yale, 1885.  
Instructor math. and civil eng, Nebraska, 1880-84, 
assoc.prof.civil eng.84-90, prof, 90-93; visited G\"ottingen and Berlin,
1898-1899; math, Stanford, 1893-1901;
civil eng, Moscow, Idaho, from 1901 , dean, col. eng, from 1911., 
died August 31, 1923 \cite{Amer,Yale}.}, made a considerable progress
on the enumeration problem so that by 1900 there were in existence tables
of (prime) knots up to ten crossings \cite{Ta,Kirk-1,Lit-0,Lit-1}. 
These tables were partially extended
in M.G.~Haseman's doctoral dissertation of 1917/8\footnote{Mary Gertrude
Haseman, born March 6, 1889, Linton, Indiana, was the fifth doctoral 
student of C.A.Scott at Brynn Mawr College. She was teaching at
University of Illinois, and died April 9, 1979.}, 
\cite{Has-1}. Knots up to 11
crossings were enumerated by John H.~Conway \cite{Co-1} before 1969
\footnote{K.A.Perko, a student of Fox at Princeton, and later a lawyer from
New York, observed a duplication in the tables: two 10-crossing diagrams
represented the same knot, see Figure 3.8. Perko corrected also the Conway's
eleven crossing tables: 4 knots were missed \cite{Per-0,Per-2}.}.

\begin{center}
\centerline{\psfig{figure=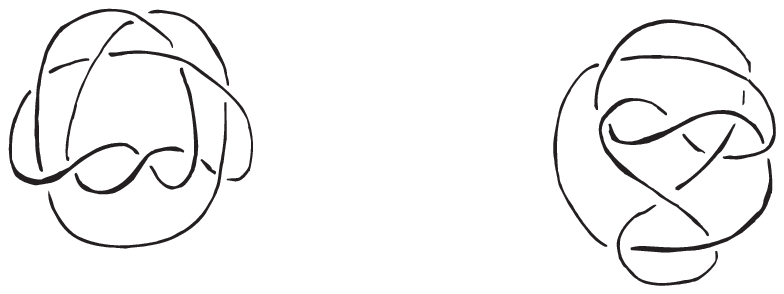,height=3.9cm}}
Figure 3.8;\  Knots of Perko.
\end{center}

Knots up to 13 crossings were enumerated by C.H.Dowker and 
M.B.Thistlethwaite \cite{D-T,This-1}, 1983\footnote
{Before my book in Polish was published in 1995, I asked Jim Hoste 
about the actual status of knot's tabulation and he 
has informed me that he and, independently, M.B.Thistlethwaite
are working on the extension up to 15 crossings of the existing knots' tables. 
They found already that there are 19 536 prime alternating, 14 crossing 
knots and 85 263 prime alternating, 15 crossing knots. Census for
all knots (not necessary alternating) is not yet verified but suggests
over 40 000 knots of 14 crossings and over 200 000 knots.}.
For the further progress we can refer to \cite{H-T} and \cite{Hos-20}. 
The number of prime, unoriented, nonalternating knots per crossing number 
$7\leq n \leq 16$ is: 0, 3, 8, 42, 185, 888, 5110, 27436, 168030, 1008906.\\
The number of prime, unoriented, alternating knots per crossing number
$3\leq n \leq 23$ is: 1, 1, 2, 3, 7, 18, 41, 123, 367, 1288, 4878, 19536, 
85263, 379799, 1769979, 8400285, 40619385, 199631989, 990623857, 4976016485, 
25182878921.\footnote{ Ortho Flint and Stuart Rankin, with coding 
by Peter de Vries, calculated alternating(23) = 25182878921 on a Compaq ES 45 
in 228 hours, finishing on Mar 14, 2004 \cite{EIS}.}
Knots and their mirror images are not counted separately.

To be able to make tables of knots,
Tait introduced three basic principles (called now the Tait conjectures).
All of them have been solved. The use of the Jones polynomial makes
the solution of the first two Tait conjectures
astonishingly easy \cite{M-4,This-3,K-6} and the
solution of the third Tait conjecture also uses essentially Jones type
polynomials \cite{M-T-1,M-T-2}.
We formulate these conjectures below:
\begin{description}
\item [T1.] 
An alternating diagram with no nugatory crossings, of an alternating
    link realizes the minimal number of crossings among all diagrams
    representing the link. A nugatory crossing is drawn (defined) in 
    Figure 10(a).
\item [T2.]
 Two alternating diagrams,   with no nugatory crossings, of the same
    oriented link have the same Tait (or writhe) number, i.e. the signed sum 
    of all crossings of the diagram with the convention that \ 
\parbox{0.6cm}{\psfig{figure=L+maly.eps}}\  is $+1$ and \ 
\parbox{0.6cm}{\psfig{figure=L-maly.eps}}\   
is $-1$.
\item [T3.]
 Two alternating diagrams, with no nugatory crossings,  of the same
    link are related by a sequence of flypes (see Figure 3.9(b)).
\end{description}


\begin{center}
\centerline{\psfig{figure=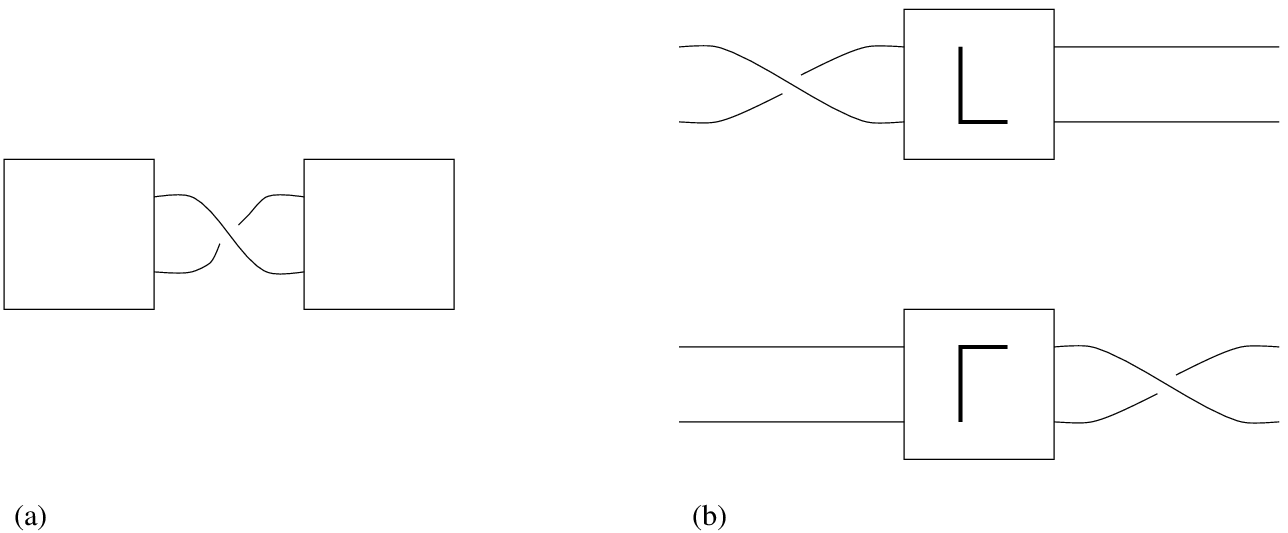}}
\end{center}

\begin{center}
                             Figure 3.9
\end{center}

A very interesting survey on the developments in knot theory in XIX
century can be found in the Dehn and Heegaard article in the Mathematical
Encyclopedia \cite{D-H}, 1907\footnote{
According to Dehn's wife, Mrs. Toni Dehn, {\it Dehn and Heegaard met at the 
third International Congress of Mathematicians at Heidelberg in 1904 
and took to each other immediately. They left Heidelberg on the same train,
Dehn going to Hamburg  and Heegaard returning to Copenhagen. They decided
on the train that an Encyclopedia article on topology would be desirable,
that they would propose themselves as authors to the editors, and that 
Heegaard would take care of literature whereas Dehn would outline a 
systematic approach which would lay the foundations of the discipline 
\cite{Mag}.}}. In this context, the papers of Oscar Simony 
from Vienna
\footnote{ Born April 23, 1852 in Vienna, died April 6, 1915 in Vienna 
\cite{Pogg}.} are of great interest \cite{Sim,Ti-2}. Figure 3.10 describes
torus knots of Simony. Simony was using continued fractions to describe 
torus knots \cite{Sim,Grg} (in essence his method was analogous to that 
employed by J.H.Conway to describe 2-bridge knots \cite{Co-1}).

\centerline{\psfig{figure=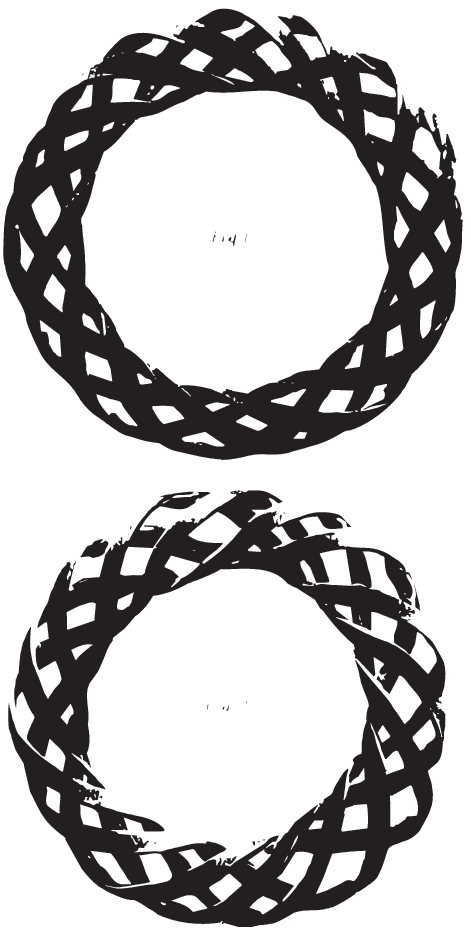,height=12.5cm}}
\begin{center}
 Figure 3.10;\ Torus knots of Simony from 1884.
\end{center}

\section{Algebraic topology in Knot Theory}
   The fundamental problem in knot theory is to be able to distinguish
non-equivalent knots. It was not achieved (even in the simple case of the
unknot and the trefoil knot) until 
Jules Henri Poincar\'e (1854-1912) in his ``Analysis Situs" 
paper (\cite{Po-1} 1895) laid foundations for algebraic topology.
Poul Heegaard (1871-1948) in his Copenhagen Dissertation
of 1898 (\cite{Heeg}) constructed
the 2-fold branch cover of a trefoil knot and showed that it is 
the lens space, $L(3,1)$, in modern terminology. He also showed that the
analogous branch cover of the unknot is $S^3$. He distinguished $L(3,1)$
from $S^3$ using the Betti numbers (more precisely he showed that the first
homology group is nontrivial and he clearly understood that it is a 3-torsion
group). 
He didn't state however that it can be used to
distinguish the trefoil knot from the unknot; see \cite
{Stil} p.226.\footnote{For the English translation of the 
topological part of the Heegaard thesis see the appendix to \cite{P-22}. }
Heinrich Tietze (1880-1964) used in 1908 the fundamental group of the
exterior of a knot in $R^3$, called the knot group, to distinguish the unknot
from the trefoil knot \cite{Ti-1}. 
The fundamental group was first\footnote
{According to \cite{Ch-M} Hurwitz' paper
of 1891 \cite{Hur} ``may very well be
interpreted as giving the first approach to the idea of a fundamental
group of a space of arbitrarily many dimensions."} introduced by
Poincar\'e in his 1895 paper \cite{Po-1}.

 Wilhelm Wirtinger (1865-1945) in his lecture 
delivered at a meeting of the German Mathematical Society in 1905 outlined
a method of finding a knot group presentation (it is called now the
Wirtinger presentation of a knot group) \cite{Wir}. Max Dehn (1878-1952),   
in his 1910
paper \cite{De-1} refined the notion of the knot group and in 1914 was able to
distinguish the right handed trefoil knot 
(\parbox{1.2cm}{\psfig{figure=+trefoilmaly.eps}}\ 
 ) from its mirror image,  
the left handed trefoil knot 
(\parbox{1.2cm}{\psfig{figure=min-trefoilmaly.eps}}\   ); 
that is Dehn showed that the trefoil
knot is not amphicheiral \cite{De-2} \footnote
{In 1978, W.Magnus wrote \cite{Mag}: 
{\it Today,  it appears to be a hopeless task to assign priorities for the
definition and the use of fundamental groups in the study of knots,
particularly since Dehn had announced \cite{De-0} one of the important 
results of his 1910 paper (the construction of Poincar\'e spaces with the 
help of knots) already in 1907.}}.

   Tait was the first to notice the relation between knots and
planar graphs. He colored the regions of the knot diagram alternately white
and black (following Listing) and constructed the graph by placing a vertex
inside each white region,   and then connecting vertices by edges going through
the crossing points of the diagram (see Figure 4.1)\cite{D-H}. 

\ \\
\centerline{\psfig{figure=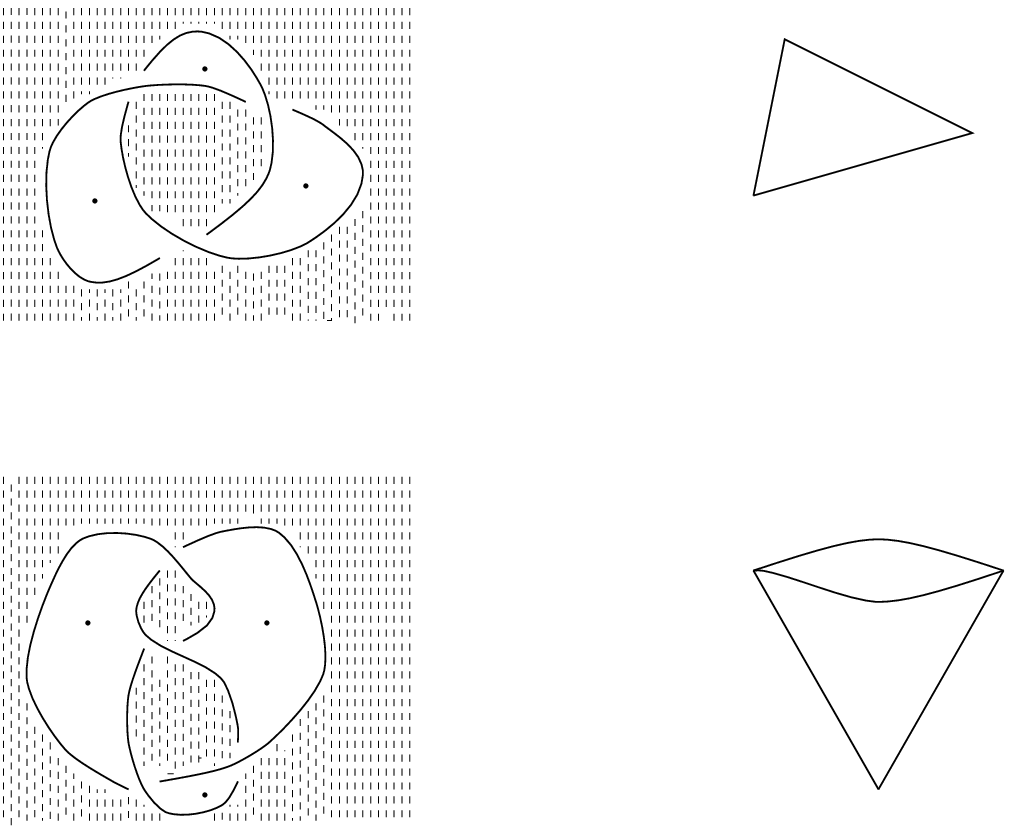,height=10.5cm}}
\begin{center}
                           Figure 4.1
\end{center}

In 1912, George David Birkhoff (1884-1944) when
trying to prove the four-color problem (formulated in 1852 by Francis
Guthrie (1831-1899)), introduced the chromatic polynomial
of a graph \cite{Birk-1}.

   The breakthrough, from the point of view which focuses on the Jones type 
link invariants, was the invention by James Waddell 
Alexander (1888-1971) of a Laurent polynomial invariant of links 
(\cite{Al-3}, 1928)\footnote{ Alexander described the numerical
precursor to his polynomial, 
for the first time, in the letter to Oswald Veblen, 1919 \cite{A-V}.}. 
Alexander was a colleague of
Birkhoff and we can conjecture that he knew about
the chromatic polynomial.\footnote{Birkhoff writes in \cite{Birk-3} :
``...Alexander, then [1911] a graduate student [at Princeton], began to be 
especially interested in the subject [analysis situs]. His well known 
``duality theorem," his contributions to the theory of knots, 
and various other results, have
made him a particularly important worker in the field". We can also mention 
that in the fall of 1909 Birkhoff became a member of the faculty at Princeton 
and left for Harvard in 1912. His 1912 paper \cite{Birk-1} ends with ``Princeton
University, May 4, 1912."}   We know for sure that 
when W.~T.~Tutte was generalizing the chromatic polynomial in 1947 \cite{Tut-1},
he was motivated by the knot polynomial of Alexander. 
The Alexander polynomial can be
derived from the group of the knot (or link). This point of view has been 
extensively developed since Alexander's discovery.
More generally, the study of the fundamental group
of a knot complement and the knot complement alone
was the main topic of research in knot theory
for the next fifty years, culminating in 1988
in the proof of
Tietze \cite{Ti-1}   conjecture (that a knot is determined by its complement) 
by Gordon and Luecke 
\cite{G-Lu}.
We can refer to the survey articles by Ralph Hartzler Fox (1913 -1973) 
\cite {F-1} and Gordon \cite{Gor} 
or books \cite{Bi-1,B-Z,K-3,Re-2,Ro-1} in this respect.
 However, Alexander observed
also that if three oriented links, $  L_+,  L_- $ and  $L_0,$  
have diagrams which are identical except near one crossing where they 
look as in Figure 4.2, then their polynomials are
linearly related \cite{Al-3}. An analogous discovery about the chromatic 
polynomial
of graphs was made by Ronald M.~Foster in 1932 (see \cite{Whit-1}; 
compare also \cite{Birk-2} Formula (10)). In early 1960's,  
 J.~Conway rediscovered Alexander's
formula and normalized the Alexander polynomial,   $\Delta _L(t) \in Z[t^{ \pm  1/2}]$,  defining it 
recursively as follows (\cite{Co-1}):
\begin{description}
\item
[(i)] $\Delta _o(t) = 1 $,   where $o$ denotes a knot isotopic 
to a simple circle

\item
[(ii)] $$\Delta_{L_+} - \Delta _{L_-} = (\sqrt t - \frac{1}{\sqrt t})\Delta _{L_0}$$
\end{description}
\ \\
\centerline{\psfig{figure=L+L-L0.eps,height=3.5cm}}
\begin{center}
Figure 4.2
\end{center}
\ \\


\section{Jones revolution}
In the spring of 1984, Vaughan Jones discovered his invariant of links, 
$V_L(t)$ \cite{Jo-0,Jo-1,Jo-2}\footnote{Jones wrote in \cite{Jo-6}:\
``It was a warm spring morning in May, 1984, when I took the uptown subway to
Columbia University to meet with Joan S. Birman, a specialist in the theory of
``braids"... In my work on von Neumann algebras, I had been astonished
to discover expressions that bore a strong resemblance to the algebraic
expression of certain topological relations among braids. I was hoping that 
the techniques I had been using would prove valuable in knot theory.
Maybe I could even deduce some new facts about the Alexander polynomial.
I went home somewhat depressed after a long day of discussions with
Birman. It did not seem that my ideas were at all relevant to the
Alexander polynomial or to anything else in knot theory.
But one night the following week I found myself sitting up in bed and running
off to do a few calculations. Success came with a much simpler approach
than the one that I had been trying. I realized I had generated a
polynomial invariant of knots."},  
and still in May of 1984 he was trying various substitutions to the 
variable $t$, in particular $t=-1$. He observed that $V_L(-1)$ is equal
to the classical knot invariant -- determinant of a knot; 
however, initially he was unable to prove it. 
Soon he realized that his polynomial satisfies
the local relation analogous to that discovered by Alexander and Conway and
established the meaning of $t=-1$.\footnote{The relation was also discovered
independently in July 1984 by Lickorish and Millett.} 
Thus the Jones polynomial is defined recursively as follows:
\begin{description}
\item
[(i)] $V_o = 1$,  
\item
[(ii)] $\frac{1}{t}V_{L_+}(t)- tV_{L_-}(t) = 
(\sqrt t - \frac{1}{\sqrt t})V_{L_0}(t)$.
\end{description}

In the summer and the fall of 1984,   the Alexander and the
Jones polynomials were generalized to the
skein (named also Conway-Jones, Flypmoth,  Homfly, Homflypt\footnote{Homfly 
or Homflypt is the acronym after the initials
of the inventors: Hoste,
Ocneanu, Millett, Freyd, Lickorish, Yetter, Przytycki and Traczyk.},
 generalized Jones,  
2-variable Jones, Jones-Conway,   Thomflyp, 
twisted Alexander) polynomial, $P_L \in Z[a^{\pm 1}, z^{ \pm 1}]$,  of oriented 
links \cite{FYHLMO,P-T-1}. This polynomial is
defined recursively as follows \cite{FYHLMO,P-T-1}:
\begin{description}
\item [(i)]
$P_o = 1$;
\item [(ii)]
$aP_{L_+} + a^{-1}P_{L_-} = zP_{L_0}.$
\end{description}

In particular $\Delta _L(t) = P_L(i,  i(\sqrt t - \frac{1}{\sqrt t}))$, 
$V_L(t) = P_L(it^{-1},  i(\sqrt t - \frac{1}{\sqrt t}))$.
In August 1985 L.~Kauffman found another approach to the Jones 
polynomial\footnote{First he thought that he produced 
a new knot polynomial and
only analyzing the polynomial he realized that he found a variant of
the Jones polynomial. }.
It starts from
an invariant,   $<D> \in Z[\mu, A, B]$,  of an unoriented link diagram 
$D$ called the Kauffman bracket \cite{K-6}.
The Kauffman bracket is defined recursively by:
\begin{description}
\item [(i)]
\[ < \underbrace{ o \ldots o}_{ i }> = \mu^{i-1} \]
\item [(ii)]
$$ <L_+> = A <L_0> + B <L_{\infty}>$$
\item [(iii)]
$$<L_-  > = B <L_0> + A <L_{\infty}> $$
\end{description}
where $L_+, L_-, L_0$ and $L_{\infty}$ denote four diagrams that are identical
except near one crossing as shown in Figure 14,   and 
\mbox{$ < \underbrace {o \ldots o}_i>$}
denotes a  diagram of $i$ trivial components ($i$ simple circles).

If we assign $B = A^{-1}$ and $\mu = -(A^2 + A^{-2})$ then the Kauffman bracket gives
a variant of the Jones polynomial for oriented links. Namely,   
for $A = t^{-\frac{1}{4}}$ and $D$ being an oriented diagram of $L$ we have
\begin{equation}
V_L(t) = (-A^3)^{-w(D)} <D>
\end{equation}
where $w(D)$ is the {\em planar writhe} ({\em twist} or {\em Tait number})
of $D$ equal to the algebraic sum of signs of crossings.

It should be noted,   as first observed by Kauffman,  
that bracket $<\ >_{\mu ,  A,  B}$
is an isotopy invariant of alternating links (and their connected sums) 
under the assumption  that the third Tait conjecture (soon after proven 
by Menasco and Thistlethwaite \cite{M-T-1,M-T-2}) holds.

In the summer of
 1985 (two weeks before discovering the ``bracket"),   L.~Kauffman invented
 another invariant of links \cite{K-5},\   
 $F_L(a,  z) \in Z[a^{ \pm  1}, z^{ \pm  1}]$,  
generalizing the polynomial discovered at the 
beginning of 1985 by Brandt,   Lickorish,  Millett and Ho
\cite{B-L-M,Ho}. To define the
 Kauffman polynomial we first introduce the polynomial
invariant of link diagrams $\Lambda _D (a, z)$.
It is defined recursively by:
\begin{description}
\item
[(i)] $\Lambda _o (a,  z) = 1$, 
\item [(ii)] $\Lambda_{{\psfig{figure=R+maly.eps}}}
 (a,  z) = a \Lambda_| (a, z); 
\ \Lambda_{{\psfig{figure=R-maly.eps}}}  
(a, z) = a^{-1} \Lambda_| (a,  z)$, 
\item [(iii)]
$\Lambda_{D_+}(a,  z) + \Lambda_{D_-}(a, z) = z(\Lambda_{_0}(a, z) + 
\Lambda_{D_\infty}(a, z))$.
\end{description}

The Kauffman polynomial of oriented links is defined by
$$F_L(a,  z) = a^{-w(D)} \Lambda _D(a, z)$$
where $D$ is any diagram of an oriented link $L$.

\ \\
\centerline{\psfig{figure=L+L-L0Linf.eps,height=2.6cm}}
\begin{center}
Figure 5.1
\end{center}

   Jones type invariants lead to invariants of three-dimensional
manifolds \cite{At,Wit,P-5,Tu-2,R-T-2,Tu-Vi,Tu-We,P-19}.
   We already mentioned that Jones type invariants of knots have been
used to solve Tait conjectures. The Jones discovery, however, not only 
introduced a delicate method of analyzing knots in 3-manifolds but
related knot theory to other disciplines of mathematics and 
theoretical physics, for example statistical mechanics, quantum field
theory, operator algebra, graph theory and computational  complexity.
On the other hand the Jones polynomial gives a simple tool to recognize
knots and as such is of great use for biologists (e.g. for analysis of
DNA) and chemists (see for example \cite{SCKSSWW}).
  \\ \ \\ 
{\bf Biographical Notes}
\\ \ \\
We give below a chronological list  of selected people mentioned in the 
article (compare \cite{B-L-W}):\\
\ \\
Gottfried Wilhelm Leibniz       (1646-1716)\\
Heinrich Kuhn                   (1690-1769)\\
Leonhard Euler                  (1707-1783)\\
Alexandre-Theophile Vandermonde (1735-1796) \\
Carl Friedrich Gauss            (1777-1855) \\
Thomas Penyngton Kirkman        (1806-1895)\\
Johann Benedict Listing         (1808-1882)\\
Hermann Von Helmholtz           (1821-1894)\\
Gustav Robert Kirchhoff         (1824-1887)\\
William Thomson (Lord Kelvin)   (1824-1907)\\
Peter Guthrie Tait              (1831-1901)\\
James Clerk Maxwell             (1831-1879)\\
Oscar Simony                    (1852-1915)\\
Jules Henri Poincar\'e          (1854-1912)\\ 
Charles Newton Little           (1858-1923)\\
Hermann Karl Brunn              (1862- 1939)\\
Wilhelm Wirtinger               (1865-1945)\\ 
Poul Heegaard                   (1871-1948)\\ 
Max Dehn                        (1878-1952)\\ 
Heinrich Tietze                 (1880-1964)\\
George David Birkhoff           (1884-1944)\\ 
James Waddell Alexander         (1888-1971)\\
Mary Gertrude Haseman           (1889-1979)\\
Kurt W.~F. Reidemeister          (1893-1971)\\
Ronald M.~Foster                (1896-1998)\\
Hidetaka Terasaka               (1904-1996)\\
Herbert Seifert                 (1907-1996)\\
Ralph Hartzler Fox              (1913-1973).\\
\baselineskip=9pt
 \small


\ \\ \ \\ \ \\
\noindent \textsc{Dept. of Mathematics, Old Main Bldg., 1922 F St. NW \\
The George Washington University, Washington, DC 20052}\\
e-mail: {\tt przytyck@gwu.edu}

\newpage
{\bf \large Appendix: Preface to P.~Heegaard thesis of 1898}\ 

Poul Heegaard (1871-1948) is well known to topologists because of his
theorem about decomposition of any $3$-manifold into two handlebodies.
His achievements in knot theory are much less known.\footnote
{The only, known to me, research papers on topology by Heegaard are
\cite{D-H,Heeg}. For other Heegaard's work see his bibliography
(Brjan Toft found the box donated to Odense Math. Dept. by
Niels Erik N{\o}rlund (1885-1981) with Heegaard papers (as in the
Heegaard's bibliography except \cite{D-H}
and only French version of \cite{Heeg}).} One of the reasons
for this was the fact that his 1898 Doctoral Dissertation \cite{Heeg}
was written in Danish. It was translated into French only 18 years later.
I learned of the Heegaard work on knots only from the book by
J.~C.~Stillwell \cite{Stil}. When visiting Odense University (1992-1994), 
I had a privilege to have every day access to a copy of Heegaard 
dissertation and I was lucky to know very gifted high school student, 
Agata Przybyszewska, who kindly agreed to translate parts of the dissertation.
We present here the Przybyszewska's translation
(from Danish) of the Preface, and the section 6 and 7
of Heegaard's Dissertation. The second, topological part of the 
dissertation is translated by Przybyszewska in [P-22] 
(references in the main part of the chapter).
One can hope that the whole Dissertation will be translated soon.
As to importance of Heegaard results we can best quote from
the Stillwell book [Stil]:\\
``Heegaard's results lay dormant (although noted by Tietze 1908) until
the publication of the French translation of his thesis in 1916.
The translation was checked for mathematical soundness by J.W.Alexander,
fresh from his work on homology groups, and we may surmise that the collision
of these ideas led to the fruitful discoveries which were to follow.
Alexander must also have read Tietze 1908 at this time, because in short
order he disposed of two of the most important of Tietze's conjectures:
Alexander 1919a\footnote{Note on two three-dimensional manifolds
with the same group, {\em Trans. Amer. Math. Soc.}, 20 (1919), 339-342.} 
shows that there are nonhomeomorphic lens spaces with the 
same group, while Alexander 1919b\footnote{Note on Riemann spaces,
{\em Bull. Amer. Math. Soc.}, 26 (1919), 370-372.}
 proves that any orientable 3-manifold
is a branched cover of $S^3$. Later in 1920 he finally took the cue
from Heegaard's example and looked for torsion in cyclic covers of $S^3$
branched over various knots."
 \\ \ \\ \ \\ \ \\ \ \\
\newpage
\ \\ 
\centerline{\bf Preliminary studies towards the topological theory of}
\centerline{\bf connectivity of algebraic surfaces.}
\centerline {Poul Heegaard}
\centerline { Copenhagen, 1898} 
\ \\
{\bf 1. Preface} \\
\ \\
It is commonly known, which development took the theory of functions of 
one independent variable, when the 
imaginary values of the independent variable were considered and the theory
 became linked with a geometrical presentation of the imaginary numbers.
 
When we shall construct a theory of functions with two independent 
variables, 
it will therefore be natural to look for a similar presentation. 
The fact that such a long time has lapsed before the work on this 
quite natural generalization was started can, 
among other reasons, be explained by the fact that the examinations with two
 independent variables are a lot more difficult than with one. 
The variety of possibilities produces conditions, for which there is no 
analogy in the theory 
of one independent variable. Picard observes thus:
 
`` On voit, par ce qui pr\'ec\`ede, les diff\'erences profondes qui separent 
la th\'eorie 
des fonctions alg\'ebrique d'une variable de la th\'eorie 
de fonctions alg\'ebriques de deux 
variables ind\'ependantes. L'analogie qui souvent est un guide excellent, 
peut devenir ici bien 
trompeuse."\footnote{We see from the proceeding the important 
differences 
that separate the theory of algebraic functions of one variable from the 
theory of 
algebraic functions of two independent variables. 
The analogy that often is an excellent guide, can just as well be faulty.}
 
Another difficulty on the way  towards a lucid  presentation 
is of course the fact that the geometrical structures that should play the 
part of the Riemann 
surfaces are 4-dimensional.

Even if the analogy can be misleading in our examination of the details, 
then a survey of the methods employed in the theory of one independent 
variable will give us a good working program for such a study. 
Let us then present such a survey.

The theory of functions with one independent variable is very closely 
connected with the theory of the algebraic curves. The geometry of such 
a curve becomes therefore of fundamental importance.

The examinations are set up from quite a different points of view. 
The most important made use of:
\\ \ \\
1. the elementary algebraic theorems. Among these we can mention
\begin{enumerate}
\item[(a)] Examination of the adjoint polynomials, created by Brill and 
           N\"other in the paper: Ueber die algebraischen Functionen 
           (M.A. vol. 7, 1873).
\item[(b)] The examinations of the linear groups of points made by Italians. 
They have tried to liberate the former theory from its projective form, 
so they could develop it independently of such concepts as degree, 
class e.t.c., cf. a summary by Castelnuovo and Euriques 
(Sur quelques recent resultantas dans la th'eorie des surfaces 
alg\'ebriques. M.A. vol. 48, p.242, 1897)
\item[(c)] Number-geometrical examinations, especially a number of papers 
by Zeuthen (e.g. M.A. vol.3 and vol. 9)
\end{enumerate}
2. the examinations of the algebraic curves which belong to 
the transcendental functions. These examinations originate from 
Riemann's pioneer work (1857) are well known, so there is no reason 
to mention them any further.\\
3. the topological examinations of the Riemann surfaces that represent 
the algebraic curve. Here once again two different approaches were chosen:
\begin{enumerate}
\item[(a)] We can either determine the connectivity number of surfaces 
using a theorem which we can regard as a generalization of Euler's 
theorem about polyhedra (Riemann, Neumann).
\item[(b)] Or we can puncture the Riemann surface and next bring it 
by continuous deformation into a normal form. As far as we can tell 
this approach was only implemented by Jul. Petersen (Foreleasninger 
over Funktionsteori, chapter IV). Listing uses indeed such a procedure 
in shaping the ``diagram" of a  spatial figure in hereby arriving to 
a generalization of Euler's theorem (Census r\"aumlicher Complexe, 1862), 
and Betti (1871) uses such a consideration in his examinations of 
the connectivity number of a n-dimensional space in generality, 
but both of these papers seem to have been unnoticed for a long time. 
This method gives us a remarkably lucid presentation of 
the discussed situation.
\end{enumerate}
The transformations of algebraic surfaces play an analogous part to 
the theory of functions of two variables. There already exist an amount 
of works -  especially from recent times - in which the problem 
is treated from points of view that correspond to these enumerated here.  \\
1. Elementary algebraic examinations.
\begin{enumerate}
\item[(a)]
Adjoint polynomials. This theory originates from Clebsch (C.R. Dec. 1868)
an N\"other (Zur aindeutigen Entsprechen... I $\&$ II, M.A. vol.2, 1869
and vol.8, 19874)
\item[(b)] 
Linear systems of curves. The Italians have created a theory of linear
systems of curves on surfaces that is analogous with the before mentioned
theory of groups of points on curves. By means of this the invariants
of surfaces can be determined (cf. the mentioned summery by Castelnuevo 
and Enriques).
\item[(c)]
Surface transformations have also been examined by Zeuthen from the number
geometrical point of view (\'Etudes g\'eom'etriques..., M.A. vol.4).
\end{enumerate}
2. Examinations by transcendent functions. Already in his paper in  
M.A. vol. 2 N\"other considers the integrals of the form:
$$ \int\int\frac{Q(xyz)dxdy}{f'_x}$$
Picard introduces the integrals of the form:
$$\int_{x_0y_0z_0}^{xyz} Pdx +Qdy$$
Where P and Q satisfy the condition of integrability (Liouville 
Journ. $1885\&1886$). Finally, Picard has in his prized paper (Memoire 
sur la the\'orie des fonctions alg\'ebriques de deux variables 
(Liouv. Journ. series 4. vol. 5, 1889) and in the recently printed 
book on the same subject (Picard at Simart: Th\'eorie des fonctions 
alg\'ebriques de deux variables ind\'ependantes, 1987)) given a coherent 
presentation of the whole theory.\\
Poincar\'e's paper "Sur les residus..." (Acta mathematica, vol. 9, 1887) 
should also be mentioned here.\\
\ \\
3. Topological examinations. This direction it is almost empty. 
In the works of Picard there is a lot, but none of it is carried through, 
as he whenever possible prefers the analytical presentation. 
The whole trouble is that the Riemann-Betti theory of connectivity numbers 
is very insufficient and difficult to use when we speak of manifolds of 
more than 2 dimensions. Poincar\'e has tried to complete it 
(Analysis Situs, Journ, de l'\'ecole polit\'echnique, series 2. Cah. 1. 1895) 
without succeeding, according to our opinion. 
Later Picard has given his share: W.~Dyck has already before worked 
on the problem (Beitri\"age sur Analysis situs, M.A. vol. 32 and 37), 
but a quite satisfying theory can not be found anywhere. \\
\ \\

Previous to the examination in this direction, we therefore need to have 
a theory of topological correspondence of manifolds of dimension greater 
than 2. What already does exist in this direction, can be compared to that 
what is mentioned in 3(b).

The following pages contain no unified totality, only a study of 
the problem; its difficulty can be used as an excuse therefore. 
To let the examinations be considered as they should be, it will perhaps 
be a good idea to advance a line of thought that have been the guide 
of my examinations.

By accident I have noticed that the Zeuthen-Halphen generalization of 
the genus theorem could be proved by a pure topological argument; 
namely, for any two Riemann surfaces, for which we assume a $\mu - \nu$ 
value-correspondence between their points, we can construct new Riemann 
surfaces, which correspond bijectively to each other. The equation, which 
says that two constructed surfaces have equal connectivity number, states 
exactly the mentioned theorem. The fact, which is so important to 
the enumerative geometry, that the generalized genus theorem can be used 
without infinitesimal examinations of the coincident points 
(or their substitutes), while these examinations are necessary 
when using the other formulas of correspondence (cf. Acta mathem., 
vol. 1, p. 171) let us see these circumstances in a new light. 
As I have found later the proof in a paper written by 
Hurwitz (M.A. vol. 39), I shall not go into details of the proof here.

My thought was that it could be possible to construct a theorem similar 
to that for the algebraic surfaces, but before that a lot of work ought 
to be done: there should be created something for algebraic surfaces 
that could correspond to the Riemann surfaces of algebraic curves; 
next topological criteria for their unique correspondence should be advanced. 
And, as it is already mentioned, the material for such an examination 
that already existed was either insufficient or full of mistakes. 
It was in trying to correct these mistakes, and to fill up the gaps, 
that the contents of the following pages originated. \\ 
\ \\ \\ 
\ \\ 
{\bf Second part}\\
On topological connectivity numbers \\
\ \\
{\bf 6 Topology.}\\
\\
Descartes' analytical method was presumably a universal method for
solving geometrical problems, but as a rule it gives constructions that
are far behind the Greeks in simplicity and elegance. During his attempts
of penetrating the principles of the curious geometrical analysis,
Leibniz formulated a number of observations, which he called
{\it Analysis situs} or {\it Geometria situs}.
Analysis situs or Geometria situs.
(Leibnizens gesammelte Werke herausg. von G. H. Pertz,1858, mathematische
Schriften vol. 1, p. 178: De analysi situs). Sometimes in connection
with this, Leibniz has been named the father of the modern topology:
i.e. the branch of mathematics, which aims at the qualitative properties
of the objects analyzed without occupying itself with the quantities and
the metrics. Leibniz notices somewhere on the occasion of his 
theory:\footnote{A figure generally contains besides an extent
(quantitas) a nature (qualitas) or form.}\\
\centerline{``Figura in universum praeter quantitatem continet qualitatem
seu formam,"} 
but in reality his theory has no points of similarity with what we
today understand as topology. In the recent times the mathematicians
became interested in a lot of topological problems. It should be                enough to remind oneself of
\begin{itemize}
\item 
the studies of topological knots, nets, and the like by Tait or Simony;
\item
the graphs;
\item
the appearance of the graphical curves; the problem of colorings,
and hereby separating, the countries on a map by means of 4 different colors;
\item
the problem of folding a stamp
\item
and ``last but not least" of the extensive attempts to make a theory
of the connectivity numbers of the n-dimensional manifolds, and of
generalizing Euler's theorem of polyhedra -- 2 attempts that are closely
connected.
\end{itemize}
\ \\
(A good survey of the literature can be found in (W.Dyck: Beitr\"age
zur Analysis situs, M.A. vol. 32)). The usage is somewhat staggering,
but it is most fair to reserve the name `analysis situs' for the
last mentioned kind of studies, and against it, to use the
name `topology' about all research of qualitative nature, as Listing
has proposed it (Vorstudien zur Topologie). \\ \ \\ \ \\
{\bf 7 Analysis situs.}\\
\ \\
Presumably, it will be very difficult to formulate basic theorems about
the connectivity in n-dimensional manifolds; anyhow it will probably be
a good idea, if we did study the concrete cases in more details,
than it has been done so far, instead of going at the whole matter
in its abstract generality. In theory, the logic should be enough to
ensure the development of mathematics, but practice shows us,
what a mighty lever is the sense of tact that develops when a theory 
is applied to
a large class of the concrete cases; even if a theory appears to be very
logical and plausible, it can very easily contain mistakes that will only
be revealed when it is additionally controlled by
the understanding.

{\it Riemann} {\it and Betti} are the first who have tried to generalize
the theory of connectivity of surfaces, which Riemann himself with great
luck has applied to the theory of Abelian integrals.
After Riemann's death Betti put
forward a paper on the mentioned subject (Sugli spazi di un numero qualanque
di dimensioni, Ann. di matem. series 2, vol. 4, 1871). He does not speak of
any collaboration between himself and Riemann; but we judge from
the fragments of the theory, which Weber has collected together from
notes that Riemann wrote down (Riemann: Gesammelte mathematische
Werke, 2. Aufl. fragment XXIX), that Riemann has contributed essentially
during his stay in Italy (Ges. Werke, p.555) to the thoughts that are
found in Betti's paper. It is where the definition can be found of
the connectivity numbers of the manifolds of dimensions greater than 2.
Later, different people have tried to improve and complete the theory:
\begin{itemize}
\item
Dyck in "Beitr"age zur Analysis situs"(M.A. vol. 32 and vol.37),
\item
Poincar\'e in the paper ``Analysis Situs", and
\item
Picard in his studies of functions of two variables 
(see p.5)
\end{itemize}

Even before I knew the last mentioned works, I have decided to try another
way than the one of Riemann-Betti that is to say, an attempt to generalize
{\it Jul. Petersen's} puncture method which I recalled from lectures.
(I became acquainted with Listing's ``Diagram" and ``Trema" a long time
after, and the same goes for the studies by Betti,
which are related to this subject, and which I until lately only have
known from the short summary in ``Fortshritteder Mathematik").
Among other things it seemed awkward to me that the
$n$ connectivity numbers were not adequate for a topological characterization
of a manifold when $n > 2$. When I became acquainted with the papers
that I'll mention later, especially that one of Poincar\'e, I started
to doubt about the rightness of my choice,
as I compared the elegant methods that I met there, with the somewhat
clumsy and awkward theory that I myself was working on;
but then I thought to have discovered that the patch, which I have chosen,
would illuminate some circumstances that would not appear clearly by means
of the other method, and because I held a method in my hands to find
sufficient criteria for the equivalence of n-dimensional manifolds, so
I decided to continue despite of the great difficulties, which I encountered.
The fact that two manifolds are equivalent does mean that they correspond
to each other, point to point, in such a way that two arbitrary points
in one of them are close to coincide when the two corresponding points
in the other manifold are. (We will later return to the connection
between Poincar\'e's concept of hom\'emorphisme, Analysis Situs $\&2$,
and the concept of equivalence that has been defined here).
The question that arises is, how we shall cut a closed manifold
(vari\'ete ferm\'ee; Poincar\'e 1.c. $\&1$) to make it simply connected.
To solve this problem we use the following procedure: the manifold is
punctured: i.e. the elementary manifold, which constitutes the
neighborhood of a point, is removed. From the puncture arises
a boundary, and it is extended by the means of continuous deformation
in such a way that we remove more and more of the given manifold.
We continue like this, until a part of the boundary encounters other parts
of itself; in such places we stop the deformation, when the distance
between the parts that are encountering has become infinitely small.
If we continue in this way, then we end up with a diagram, which is made
from a system of manifolds of lower dimensions than the original
one\footnote{It is called a spine of the manifold in the modern
terminology [translator's remark]} -- or rather: it is made from 
the closest neighborhood of this system, that is, 
a manifold that is infinitely small in the $n^{th}$
dimension. We call the system of manifolds of 
a lower dimension, which is the boundary of the diagram, its nucleus.
This diagram has two meanings. First, it presents the manifold which is
equivalent with the given one after the puncture. If we succeed in
proposing normal forms into which the diagrams could be reduced,
then the fact that the normal forms of the diagrams are
identical will be the necessary and sufficient condition for two diagrams
to be equivalent. -- But the nucleus of the diagram indicates the cuts
that shall be done to make the given manifold simply connected;
if we make the discussed extending in reverse order, then the manifold,
which remains after performing the cuts determined by
the diagram, will contract into that elementary manifold, which is removed
by the puncture.

We will later return to this comparison of results that
are found along this path, and the Riemann-Betti theory of
connectivity numbers.\\
\ \\
The whole {\it Second part} is translated in [P-22].
\ \ \\ \ \\
{\bf Poul Heegaard's bibliography:}

[Heeg]
P. Heegaard, Forstudier til en Topologisk Teori for de
algebraiske Fladers Sammenh{\ae}ng, K{\o}benhavn, 1898,
Filosofiske Doktorgrad;
French translation: Sur l'Analysis situs, {\em  Soc. Math. France Bull.},
44 (1916), 161-242.

[D-H]
M. Dehn, P. Heegaard, Analysid situs, {\em Encykl. Math. Wiss.}, vol. III
AB3 Leipzig, 1907, 153-220.

[H-1]
P. Heegaard, Der Mathematikunterricht in D\"anemark, {\em Internationale
Matchematische Unterrichtskommission I.M.U.K}, K{\o}benhavn 1912. (Teaching
of Mathematics in Denmark).

[H-2]
P. Heegaard, Bidrag til Grafernes Theori,
Tredie Skandinaviske Matematikerkomgres i Kristiania, 1913, 1-6.

[H-3]
P. Heegaard, Omriss av kronologiens historie i Europa (med s{\oe}rlig
henblikk p\'a den norske almanakk av 1644, ????? 51-96.

[H-4]
P. Heegaard, Hamiltons Dodekaederspil, {\em Nyt Tidsskrift Matematik},
29(4), 1919, 81-88.

[H-5]
P. Heegaard, Mindetale, oer prof..dr.H.G.Zeuthen, S{\ae}rtrykk av
Videnskapsselskapets Forhandlinger 1920, 1-3.

[H-6]
P. Heegaard, Hovedlinjer i H.G.Zeuthens videnskabelige production,
Foredrag i Norsk Matematisk Forening, {\em S{\ae}rtrykk av Norsk Matematisk
Tidsskrift}, 1920, 1-5.

[H-7]
P. Heegaard, Eudoksos' Hippopede. Et foredrag holdt i Norsk Matematisk Forening,
{\em S{\ae}rtrykk av Norsk Matematisk Tidsskrift}, 1921, 1-10.

[H-8]
P. Heegaard, Felix Klein, Et foredrag holdt i Norsk Matematisk
Forening 25. mars 1926, {\em S{\ae}rtrykk av Norsk Matematisk
Tidsskrift}, 1926, 1-7.

[H-9]
P. Heegaard, Om en generalisasjon til rummet av Sophus Lies
fremstilling av imagin{\ae}re elementer i planet, Den syvende Skandinaviske
Matematikerkongress, I Oslo 19-22 August 1929, Oslo 1930, 92-95.

[H-10]
P. Heegaard, Scandinavie. Les modifications essentielles de
l'enseignement math\'ematique dans les principaux pays depuis 1910,
{\em L'Enseignement math\'ematique}, XXIX(4-6), 1930, 307-314.

[H-11]
P. Heegaard, Appendix P.Osl 6. The date of the Horoscope of Philoe (in ``Papyri
Osloenses'' by S.Eitrem, L.Amundsen), Oslo, on commission by Jacob Dybwad, 1931,
146-151.

[H-12]
P. Heegaard, \"Uber die Heawoodschen Kongruenzen, {\em S{\ae}rtrykk av Norsk
Matematisk, Forenings Skrifter},Serie II, Nr. 1-12, 1933, 47-54.

[H-13]
P. Heegaard, Norv\`ege. La pr\'eparation th\'eorique et pratique
des professeurs de Math\'ematiques de l'enseignement secondaire,
{\em De L'Enseignement math\'ematique}, 5-6, 1933, 360-364.

[E-H] F. Engel, P. Heegaard, Sophus Lies Samlede Avhandlinger. Blind I,
{\em S{\ae}rtrykk av Norsk Matematisk Tidsskrift}, 1934, 1-9.

[H-14]
P. Heegaard, Et brev fra Abel til Degen, {\em S{\ae}rtrykk av Norsk Matematisk
Tidsskrift}, 2. hefte 1935, 1-6.

[H-15]
P. Heegaard, Bemerkungen zum Vierfarbenproblem, {\em Mat. Sb.} (New Ser.),
1(43):5, (1936), 685-693.

[H-16]
P. Heegaard, La repr\'esentation des points imaginaires de Sophus Lie at sa
valeur didactique, Atti del Congresso Internazionale dei Matematici,
Bologna, 3-10 settembre 1938 - VI, 421-423.

[H-17]
P. Heegaard, Die Topologie unddie Theorie der algebraischen
Funktionen mit zwei komplexen Variabeln, Congr\`es des Math\'ematiques
\`a Helsingfors 1938, 1-8.

 \end{document}